\documentclass[11pt]{article}
\usepackage{latexsym}
\usepackage{amsfonts}
\usepackage{mathrsfs,amsthm,amsmath}
\usepackage{color}
\newtheorem {theorem}{Theorem}[section]

\newtheorem {corollary}{Corollary}[section]

\newtheorem {definition}{Definition}[section]
\newtheorem {example}{Example}[section]

\newtheorem {lemma}{Lemma}[section]

\newtheorem {proposition}{Proposition}[section]

\begin{document}
\title{An inverse Sanov theorem for exponential families\thanks{We acknowledge the support of
Indam-GNAMPA (for the Research Project \lq\lq Stime asintotiche: principi di invarianza e grandi deviazioni\rq\rq),
MIUR (for the Excellence Department Project awarded to the Department of Mathematics, University of Rome Tor Vergata (CUP E83C18000100006)),
University of Rome Tor Vergata (for the Research Program \lq\lq Beyond Borders\rq\rq, 
Project \lq\lq Asymptotic Methods in Probability\rq\rq (CUP E89C20000680005)),
and Sapienza Universit\`a di Roma (for the Project \lq\lq Modelli stocastici con applicazioni nelle scienze e nell'ingegneria\rq\rq (CUP RM11715C7D9F7762)).}}
\author{Claudio Macci\thanks{Address: Dipartimento di Matematica, Universit\`a di Roma Tor Vergata, 
Via della Ricerca Scientifica, I-00133 Rome, Italy. e-mail: \texttt{macci@mat.uniroma2.it}}
\and Mauro Piccioni\thanks{Address: Dipartimento di Matematica, Sapienza Universit\`a di Roma,
Piazzale Aldo Moro 5, I-00185 Rome, Italy. e-mail: \texttt{mauro.piccioni@uniroma1.it}}}
\maketitle
\begin{abstract}
We prove the large deviation principle (LDP) for posterior distributions arising from subfamilies of 
full exponential families, allowing misspecification of the model. Moreover, motivated by the so-called
inverse Sanov Theorem (see e.g. Ganesh and O'Connell 1999 and 2000), we prove the LDP for the 
corresponding maximum likelihood estimator, and we study the relationship between rate functions.
In our setting, even in the non misspecified case, it is not true in general that the rate functions
for posterior distributions and for maximum likelihood estimators are Kullback-Leibler divergences 
with exchanged arguments.
\end{abstract}

\section{Introduction}

The theory of large deviations deals with the asymptotics of small probabilities at 
an exponential scale. The general framework (see e.g. Dembo and Zeitouni, 1998) concerns a family of probability 
measures $\{\phi_n, n \in \mathbb{N}\}$ on some topological space $\mathcal{X}$, together with a lower semi-continuous rate function $I:\mathcal{X}\to [0,+\infty]$. The sequence $\{\phi_n, n \in \mathbb{N}\}$ satisfies the Large Deviation 
Principle (LDP) with rate function $I$ if
$$
-\inf_{x\in int(A)}I(x) \leq \liminf_{n \to \infty} n^{-1} \ln \phi_{n}(A)\leq 
\limsup_{n \to \infty} n^{-1} \ln \phi_{n}(A) \leq - \inf_{x \in {\bar A}}I(x),
$$
where $ int(A)$ and $\bar A$ are the interior and the closure of $A$, respectively. Throughout this paper the rate functions will always be good, with
compact sub-level sets $\{x\in\mathcal{X}, I(x)\leq\eta\}$, for any $\eta\geq 0$. 

In the standard case in which the function $I$ has a unique zero $x_0$, the LDP implies the convergence in probability of $\{\phi_n \}$ to the Dirac law concentrated on $x_0$. LDP's of statistical interest concern the laws of a sequence $\{T_n\}$ of consistent estimators of some parameter $\theta$. If the true value $\theta_0$ is the unique zero of $I$, choosing $A_{\varepsilon}=\{\theta: d(\theta,\theta_0)> \varepsilon\}$, $d$ being any metric compatible with the topology of $\mathcal {X}$, the above estimates yield the \emph{inaccuracy rate} of the estimator at $\theta_0$ and $\varepsilon>0$, as defined in Bahadur et al (1980) (under the regularity assumption that the infimum of $I$ on $A_{\varepsilon}$ coincide with that on its closure).

Since the rate functions presented in this paper are always expressed in terms of Kullback-Leibler (KL) divergences, we recall that the KL divergence of a probability $P$ with respect to another probability $Q$, both on the same measurable space, is defined by 
$$
D( P|| Q)=\left\{
\begin{array}{ll}
\int \log (\frac {dP}{dQ}) dP,&\ \mbox{if}\ P \ll Q,\\
+\infty,&\ \mbox{otherwise}.
\end{array}\right.
$$
The fundamental property of the KL divergence is that it is always positive, except that it vanishes when $P=Q$. A peculiar feature is that, despite $D(P||Q)$ measures the discrepancy between $P$ and $Q$, in general $D(P||Q) \neq D(Q||P)$.

The simplest form of Sanov theorem (see e.g. Theorem 2.1.10 in Dembo and Zeitouni,1998) concerns the empirical distributions $E_n=\frac {1}{n} \sum_{i=1}^n\delta_{X_i}$ when
$\{X_i, i=1,\ldots,n\}$ are i.i.d. random variables taking values on a finite set, say $\{1,\ldots, N\}$, with $P(X_1=k)=p_0(k), k=1,\ldots,N$. The Sanov theorem states that, if $\phi_n$ is the law of $E_n$
over the unit simplex $\mathcal {X}=\mathcal{X}_N$ of $\mathbb {R}^N$, then, as $n \to \infty$, $\{\phi_n, n \in \mathbb{N}\}$  satisfies the LDP principle with rate function $I(p)=D(p||p_0)$, for $p \in \mathcal{X}_N$. 
Since the empirical distribution is a \emph{Maximum Likelihood Estimator} (MLE) of the true distribution $p_0$, the Sanov theorem implies the consistency of this estimator and provides exponential rates of convergence to zero for values of the MLE different from $p_0$.

Next consider a Bayesian approach to the above problem, specifying on the space $\mathcal{X}_N$
a \emph{prior} probability measure $\nu$, supported for the sake of simplicity by the whole space. By means of the Bayes' theorem, assuming $\{X_i, i=1,\ldots,n\}$ i.i.d. conditionally to a random law drawn from $\nu$, the \emph{posterior distribution} is readily obtained: call it $\pi_N=\pi^N(\cdot|x_1,\ldots,x_n)$, where $x_i$ is the observed values of the random variable $X_i$, for $i=1,\ldots,n$. Assuming that the empirical frequencies converge to some $p_0(k), k=1,\ldots,N$, Ganesh and O'Connell (1999) proved the LDP for $\pi^N$ with rate function $J(p)=D(p_0||p)$, for $p \in \mathcal{X}_N$. This result implies \emph{Bayesian consistency from a frequentist point of view}: if the observed values are i.i.d. draws from \emph{the true law} $p_0$, the assumption on the empirical frequencies holds almost surely, and the Bayesian posterior concentrates around $p_0$ at an exponential speed, not depending on the prior $\nu$. As a matter of fact, Ganesh and O'Connell (1999) allowed a prior distribution $\nu$ whose support $S(\nu)$ strictly included in the unit simplex, assuming however $p_0\in S(\nu)$: our main result allows to generalize Ganesh and O'Connell's one when $p_0\notin S(\nu)$, i.e. the prior is \emph{misspecified}.

It is also worth to recall that the Sanov and the inverse Sanov theorem are proved in the literature under more general hypotheses. The Sanov theorem for random variables taking values on a Polish space can be found in Theorem 6.2.10 in Dembo and Zeitouni (1998). Ganesh and O'Connell (2000) proved the inverse Sanov theorem for random variables taking values in a compact space with a Dirichlet prior (see Ferguson, 1973), with a mean measure supported by the whole space. It is not surprising that in a nonparametric framework Bayesian consistency depends on the choice of the prior, as showed by many results on this topic (see the review papers of Ghosal, Ghosh and Ramamoorthi (1998) and Wasserman (1998) and Kleijn and Van der Vaart (2006), where misspecification is also discussed). Again the two LDP's have rate functions which are KL divergences with exchanged arguments.

In this paper we concentrate our attention to more general parametric models $\{P_{\boldsymbol {\theta}}, \boldsymbol {\theta} \in \Theta \}$, with $\Theta \subset \mathbb {R}^d$. In general it should not be expected that, even if a LDP holds for the corresponding MLE, the rate functions is the KL divergence $I(\boldsymbol {\theta})=D(P_{\boldsymbol {\theta}}||P_{\boldsymbol {\theta}_0})$, that will be written from now as $D(\boldsymbol {\theta}||\boldsymbol {\theta}_0)$, where $\boldsymbol {\theta}_0$ is the true value. Indeed, this would imply that the MLE is a consistent estimator that it is optimal \emph{with respect to the inaccuracy rate}. But this is well-known to be false in general: see the theoretical results in Kester and Kallenberg (1986) and Arcones (2004), and specific examples in  Macci and Petrella (2009) and Macci (2014). In the sequel we will try to give some geometrical justification for the phenomenon, in the setting that we are about to describe.

The main aim of this paper is to prove a \emph{"parametric" inverse Sanov theorem} for an exponential family $\{P_{\boldsymbol {\theta}}, \boldsymbol {\theta} \in \Theta \}$. To explain what this means recall that, once a prior distribution $\nu$ supported by $\Theta$ is specified, the sequence of regular posterior distributions is constructed according to the Bayes' theorem, assuming that a sample of size $n$ is drawn conditionally from $P_{\boldsymbol {\theta}}$, $\boldsymbol {\theta}$ being a $\nu$-distributed random variable. Then, provided the sequence of observed sample values satisfies some natural condition
(in the sense that it holds with probability $1$ under repeated i.i.d. sampling from $P_{\boldsymbol {\theta}_0}$), we show that the sequence of posterior distributions satisfies a LDP with rate function $J(\boldsymbol {\theta})=D(\boldsymbol {\theta}||\boldsymbol {\theta})$, for $\boldsymbol {\theta} \in \Theta$.

The above framework allows to discuss also some type of misspecification. Indeed any exponential family is embedded in a \emph{full exponential family}, which can be considered as a \emph{saturated model} $\{P_{\boldsymbol {\theta}}, \boldsymbol {\theta} \in \Theta_0 \}$, with $\Theta \subset \Theta_0$. When $\boldsymbol {\theta}_0 \in \Theta_0 \setminus \Theta$, we say that the model is \emph{misspecified}. In this case we still prove the LDP on $\Theta$, but the rate function becomes $J(\boldsymbol {\theta})=D(\boldsymbol {\theta}_0||\boldsymbol {\theta})-D(\boldsymbol {\theta}_0||\Pi_{\Theta}(\boldsymbol {\theta}_0))$, $\Pi_{\Theta}(\boldsymbol {\theta}_0)$ being any minimizer of $D(\boldsymbol {\theta}_0||\boldsymbol {\theta})$, with $\boldsymbol {\theta} \in \Theta$ (this is an abuse of notation, since the minimizer is not unique in general). With a terminology borrowed from information geometry (Nielsen, 2018) we say that $\Pi_{\Theta}(\boldsymbol {\theta}_0)$ is a $m$-projection of $\boldsymbol {\theta}_0$ on $\Theta$. A question we will be interested in is the relation of $J(\boldsymbol {\theta})$ with $D(\Pi_{\Theta}(\boldsymbol {\theta}_0)||\boldsymbol {\theta})$.

Once the above framework is established, the problems with the existence of a "parametric Sanov theorem" are better understood. It is well-known that, under rather general conditions, the sequence of sample means of i.i.d. random vectors satisfies a LDP, known as the Cram\'er's theorem (see e.g. Theorem 2.2.30 in Dembo and Zeitouni (1998)). Since for exponential families the sufficient statistics are sample means, as long as the MLE within the family $\{P_{\boldsymbol {\theta}}, \boldsymbol {\theta} \in \Theta \}$ is a continuous function of the sample mean, the \emph{contraction principle} (see Theorem 4.2.1 in Dembo and Zeitouni, 1998) allows to transfer the LDP from the domain to the image of the continuous function. 
In this admittedly restrictive situation we obtain a LDP for the MLE with the rate function
$$I(\boldsymbol{\theta})=\inf_{\boldsymbol {\theta}^{\prime} \in T_{\boldsymbol {\theta}}} D(\boldsymbol {\theta}^{\prime}||\boldsymbol {\theta}_{0}),$$
where $T_{\boldsymbol {\theta}}(\Theta)\subset \Theta_0$ is a certain subfamily which intersects \emph {"transversally"} $\Theta$ at $\boldsymbol {\theta}$. In information geometry, see e.g. Nielsen (2018), any argument $\boldsymbol {\theta}^{\prime}\in T_{\boldsymbol {\theta}}(\Theta)$ achieving the infimum in the above display is called 
an $e$-projection of $\boldsymbol {\theta}_{0}$ in $T_{\boldsymbol {\theta}}(\Theta)$, that will be denoted by $\Pi^{T_{\boldsymbol {\theta}}(\Theta)}(\boldsymbol{\theta}_0)$ in the sequel (a different terminology is used by Csiszar (1975), who coined the name $I$-projection instead). In the well-specified case, that is when $\boldsymbol {\theta}_0 \in \Theta$, the "parametric Sanov theorem" holds as long as this $e$-projection is actually $\boldsymbol{\theta}$. This will be investigated in the sequel and some examples will be discussed.

The plan of the paper is the following. In Section \ref{sec:preliminary} the basic tools about exponential family, convex conjugate functions and information geometry will be introduced, and we prove the main result of the paper. Moreover the consequences of misspecification are investigated. The result about the LDP for MLE estimators is presented in Section \ref{sec:LDP-MLE}. We also give a geometrical interpretation for the cases in which a "parametric" Sanov theorem holds or does not hold. 
Section \ref{sec:Letac}
deals with the notion of \emph{dual measures}, a concept arising quite naturally from the main subject of the paper (see Letac, 2022).
Section \ref{sec:dom(K)-not-open} is devoted to some extension of the main theorem for non regular families.

\section{The main result}\label{sec:preliminary}

Let us recall some general facts about natural exponential families, for which suitable references are the books of Barndorff-Nielsen (1978) and Brown (1986). Another convenient reference on this subject is the paper by Morris (1982) which however concerns only the one-dimensional case. We assume that a $\sigma$-finite Borel measure on $\mathbb{R}^d$ is specified, called $\lambda$ from now on. We assume that $\lambda$
is not concentrated on a proper linear submanifold of $\mathbb{R}^d$, and its \emph {cumulant generating function}
$$\kappa (\boldsymbol {\theta})=\log \int_{\mathbb {R}^d} e^{\boldsymbol {\theta}\cdot \boldsymbol {x}}
\lambda(d\boldsymbol {x}), \,\,\, \boldsymbol {\theta} \in \mathbb {R}^d,$$
is \emph{regular}, that is it has an \emph{open} domain of finiteness $\Theta_0 (\lambda)=\Theta_0$ (the \emph{essential domain} of $\kappa$). 
The \emph{full natural exponential family generated by} $\lambda$ is defined by
\begin{equation}\label{family}
	\frac {dP_{\boldsymbol {\theta}}}{d\lambda}(\boldsymbol{x})=e^{\boldsymbol{\theta } \cdot \boldsymbol{x}-\kappa (\boldsymbol{\theta})},\,\,\, \boldsymbol{\theta} \in \Theta_0
\end{equation}
and it is itself called regular in what follows.
It is well-known (Barndorff-Nielsen, Chapter 7) that the function $\kappa$ is strictly convex and smooth on $\Theta_0$, with gradient
\begin{equation}\label{grad}
\nabla {\kappa}(\boldsymbol{\theta})=\int \boldsymbol{x}P_{\boldsymbol {\theta}}(d\boldsymbol{x}),\,\,\, \boldsymbol{\theta} \in \Theta_0.
\end{equation}
This allows to compute the KL divergence 
\begin{equation}\label{exponentialKL}
D(\boldsymbol {\theta}_0||\boldsymbol {\theta})=
\int \log \left(\frac {dP_{\boldsymbol {\theta}_0}}{dP_{\boldsymbol {\theta}}}\right) dP_{\boldsymbol {\theta}_0}=\kappa (\boldsymbol {\theta})-\kappa(\boldsymbol {\theta}_0)-(\boldsymbol {\theta}-\boldsymbol {\theta}_0) \cdot \nabla \kappa (\boldsymbol {\theta}_0).
\end{equation}
Notice that this is the difference, computed in $\boldsymbol {\theta}$, between the function $\kappa$ and the linear approxmation given by the tangent hyperplane at $\boldsymbol {\theta}_0$.

The gradient $\nabla \kappa$ establishes a diffeomorphism of $\Theta_0$ onto its image $M_0(\lambda)=M_0$, which is increasing, that is $(\boldsymbol {\theta}-\boldsymbol {\theta}_0) \cdot (\nabla \kappa (\boldsymbol {\theta})-\nabla \kappa (\boldsymbol {\theta}_0))>0$ for any $\boldsymbol {\theta} \neq \boldsymbol {\theta}_0$. Under our regularity assumption the open set $M_0$ is the interior of the \emph{convex support} $C(\lambda)$ of $\lambda$, that is the smallest closed convex set whose complement is $\lambda$-neglectable (Barndorff-Nielsen, Theorem 9.2 and Corollary 9.6). Clearly this is the set where the sample means from any element of the exponential family take values. Thus the exponential family can be also parametrized by its mean, taking values in $M_0$. The interplay of the \emph{natural} and the \emph{mean value parametrization} is crucial for the understanding of some of the key properties of exponential families; however in what follows we have preferred to refer to the natural parametrization as much as possible.

Next, for $\boldsymbol{t} \in \mathbb {R}^d$, let us define the \emph{log-likelihood function}
$$
\ell(\boldsymbol{\theta}; \boldsymbol {t})=\boldsymbol{\theta }\cdot \boldsymbol{t} -\kappa (\boldsymbol{\theta }),\ \boldsymbol{\theta } \in \mathbb {R}^d.
$$
The terminology is due to the fact that when $\boldsymbol{t}$ is replaced by $\bar{\boldsymbol{x}}_{n}=\frac{1}{n}\sum_{i=1}^n \boldsymbol{x}_{i}$, the mean of an observed sample $\boldsymbol{x}_{1} ,\boldsymbol{x}_{2} ,\boldsymbol{\ldots } ,\boldsymbol{x}_{n}$, then $\ell(\boldsymbol{\theta};\bar{\boldsymbol{x}}_{n})$ is the (normalized) \emph{log-likelihood function} for the family \eqref{family}. For $\boldsymbol{t} \in M_0$, the function $\ell (\cdot;\boldsymbol {t})$ is bounded from above, and it has either empty or compact super-level sets $\{\boldsymbol{\theta} \in \mathbb {R}^d: \ell(\boldsymbol{\theta};{\boldsymbol{t}}) \geq s\}$ for $s \in \mathbb {R}$, all contained in $\Theta_0$ (Barndorff-Nielsen, Section 9.3). The \emph{convex conjugate} function of $\kappa$ is defined by
\begin{equation}\label{supachie}
\kappa^{\ast} (\boldsymbol{t}) = \sup _{\boldsymbol{\theta} \in \mathbb {R}^d} \ell(\boldsymbol{\theta};\boldsymbol{t})=\sup _{\boldsymbol{\theta} \in \Theta_0} \ell(\boldsymbol{\theta};\boldsymbol{t}),\,\,\,\, \boldsymbol {t} \in \mathbb {R}^d,
\end{equation}
which is a convex lower semi-continuous function. As already said $\kappa^{\ast}$ is finite on $M_0$, which coincides with the interior of its essential domain (this may contain points in the boundary, but in any case $\kappa^{\ast}$ is steep, see the next section for this definition). It is differentiable in $M_0$, and, under our regularity assumption,
the gradient $\nabla\kappa^{\ast}$ is the inverse mapping to $\nabla \kappa$. Moreover, for any $\boldsymbol{t}, \boldsymbol {\theta} \in \mathbb {R}^d$
$$\kappa^{\ast} (\boldsymbol{t})+\kappa (\boldsymbol {\theta}) \geq \boldsymbol{\theta } \cdot \boldsymbol{t},$$
and the equality holds if and only if $\boldsymbol{t}$ and $\boldsymbol{\theta}$ are \emph{a conjugate pair}, meaning that
$\boldsymbol{t} \in M_0$ and $\boldsymbol{\theta}=\nabla \kappa^{\ast}(\boldsymbol{t})$. As a consequence $\nabla \kappa^{\ast} (\boldsymbol {t})$ is the \emph{unconstrained} MLE for the parameter $\boldsymbol {\theta}$ and
\begin{equation}\label{achieve}
\kappa^{\ast}(\boldsymbol{t})=\nabla \kappa^{\ast}(\boldsymbol {t}) \cdot \boldsymbol {t} - \kappa (\nabla \kappa^{\ast}(\boldsymbol {t}))=\ell(\nabla \kappa^{\ast}(\boldsymbol {t});\boldsymbol {t}),\,\,\, \boldsymbol{t} \in M_0.
\end{equation}

Now let $B \subset \mathbb{R}^d$ and define
\begin{equation}\label{maximum}
\kappa _{B}^{\ast} (\boldsymbol{t}) =\sup_{\boldsymbol{\theta } \in B}\ell(\boldsymbol{\theta};\boldsymbol{t})
=\sup _{\boldsymbol{\theta } \in \mathbb{R}^d}\{\boldsymbol{\theta } \cdot \boldsymbol{t} -\kappa (\boldsymbol{\theta}) -\delta (\boldsymbol{\theta }\vert B) \},
\end{equation}
where $\delta (\boldsymbol{\theta }\vert B)=0$ if $\boldsymbol{\theta }\in B$ and $\delta (\boldsymbol{\theta }\vert B)=+\infty$ if $\boldsymbol{\theta }\notin B$. The function
$\kappa _{B}^{\ast}$ is again a lower semi-continuous convex function, being a supremum of affine functions. Given that $\kappa _{B} ^{\ast} (\boldsymbol{t}) \leq \kappa  ^{\ast} (\boldsymbol{t})$ (the former being a supremum constrained to a smaller domain), it is finite on $M_0$.
Moreover, since a convex function is continuous in the interior of its domain of finiteness (see e.g. Roberts 
and Varberg (1973), Theorem D, page 93), the function $\kappa _{B}^{ \ast} $ is always continuous in $M_0$, irrespectively of the choice of the set $B$. 
Moreover, when $B$ is closed the supremum in \eqref{maximum} is achieved in at least one point.

Next, suppose that a prior distribution $\nu$ is specified on $\Theta_0$, whose support (in the relative topology of $\Theta_0$) is $\Theta$. Since the function $\kappa$ has compact super-level sets and $\Theta$ is relatively closed in $\Theta_0$, when $\boldsymbol {t} \in M_0$ and $B=\Theta$ the supremum in \eqref{maximum}  is achieved in at least one point. Given that there exists a unique $\boldsymbol {\theta}_0$ such that $\nabla \kappa (\boldsymbol {\theta}_0)=\boldsymbol {\mu}_0$, from the expression \eqref{achieve} we have
\begin{equation}\label{remember}
D(\boldsymbol {\theta}_0||\boldsymbol {\theta})=\kappa^{\ast}(\boldsymbol {\mu}_0) - \ell (\boldsymbol {\theta}; \boldsymbol {\mu}_0),
\end{equation}
and this shows that a maximizer of $\ell (\boldsymbol {\theta}; \boldsymbol {\mu}_0)$ in $\boldsymbol {\theta} \in \Theta$ is also a $m$-projection $\Pi_{\Theta}(\boldsymbol {\theta}_0)$ of $\boldsymbol {\theta}_0$ on $\Theta$.

Now consider the following sequence of probability measures on $\Theta$
\begin{equation}
\label{bayes}
\pi _{n}(A|\bar {\boldsymbol{x}}_{n}) =\frac{\int _{A}\exp \left (n\ell(\boldsymbol{\theta };\bar{\boldsymbol{x}}_{n})\right )\nu (d\boldsymbol{\theta })}{\int _{\Theta}\exp (n\ell(\boldsymbol{\theta };\bar{\boldsymbol{x}}_{n}))\nu (d\boldsymbol{\theta })},
\end{equation}
$A$ being any Borel subset of $\Theta$, which is well defined as long as $\bar{\boldsymbol{x}}_{n} \in M_0$. Here is our main result:

\begin{theorem}\label{invsanofi} If the sequence  $\{\bar{\boldsymbol{x}}_{n}, n\in \mathbb{N}\}$ converges to some $\boldsymbol{\mu}_0 \in M_0$, with $\boldsymbol {\mu}_0=\nabla \kappa (\boldsymbol {\theta}_0)$, the sequence of probability measures $\{\pi_{n}(\cdot|\bar {\boldsymbol{x}}_{n}), n \in \mathbb {N}\}$ supported by $\Theta$ satisfies a LDP with rate function
\begin{equation}\label{invdiv}
J(\boldsymbol{\theta})=D(\boldsymbol{\theta}_0||\boldsymbol{\theta})-D(\boldsymbol{\theta}_0||\Pi_{\Theta}(\boldsymbol {\theta}_0)),\,\,\, \boldsymbol{\theta} \in \Theta,
\end{equation}
where $\Pi_{\Theta}(\boldsymbol {\theta}_0)$ is any $m$-projection of $\boldsymbol{\theta}_0$ on $\Theta$.
\end{theorem}

In order to prove this theorem, we need first to prove the following lemma.

\begin{lemma}\label{th1} Under the assumptions of Theorem \ref{invsanofi},
$$\lim_{n \to \infty} \frac{1}{n}\ln \int_{\Theta} \exp (n\ell(\boldsymbol{\theta };\bar{\boldsymbol{x}}_{n}))\nu (d\boldsymbol{\theta })
=\kappa _{\Theta}^{\ast} (\boldsymbol{\mu}_0)=\ell(\Pi_{\Theta}(\boldsymbol {\theta}_0); \boldsymbol {\mu}_0).$$
\end{lemma}
\begin{proof}
Replacing the integrand with its supremum (with respect to $\boldsymbol{\theta}$) over $\Theta$, which is the support of $\nu$, we immediately have
$$\frac{1}{n}\ln \int_{\Theta} \exp (n\ell(\boldsymbol{\theta };\bar{\boldsymbol{x}}_{n}))\nu (d\boldsymbol{\theta })\leq 
\sup_{\boldsymbol {\theta} \in \Theta} \ell(\boldsymbol {\theta};\bar {\boldsymbol {x}}_n)=\kappa _{\Theta}^{ \ast} (\bar{\boldsymbol{x}}_{n}).$$
Hence if $\bar{\boldsymbol{x}}_{n}$ tends to $\boldsymbol{\mu}_0 \in M_0$ as $n \to \infty $, then $\bar{\boldsymbol{x}}_{n}$ is eventually in $M_0$. By continuity of $\kappa_{\Theta}^{\ast}$ within this set,
$\kappa _{\Theta}^{ \ast} (\bar{\boldsymbol{x}}_{n})$ tends to $\kappa _{\Theta}^{\ast} (\boldsymbol{\mu}_0)$ and
$$\limsup _{n \rightarrow \infty }\frac{1}{n}\ln \int_{\Theta} \exp (n\ell(\boldsymbol{\theta };\bar{\boldsymbol{x}}_{n}))\nu (d\boldsymbol{\theta}) \leq \kappa _{\Theta}^{\ast} (\boldsymbol{\mu}_0).$$
For the reverse inequality, recall that
\begin{equation}\label{eq:aggiunta}
\kappa _{\Theta}^{\ast} (\boldsymbol{\mu}_0)=\ell(\Pi_{\Theta}(\boldsymbol {\theta}_0); \boldsymbol{\mu}_0),
\end{equation}
for some $\Pi_{\Theta}(\boldsymbol {\theta}_0) \in \Theta$. Suppose first that $\nu(\{\Pi_{\Theta}(\boldsymbol {\theta}_0)\})>0$. In this case 
$$
\ell(\Pi_{\Theta}(\boldsymbol {\theta}_0);\bar{\boldsymbol{x}}_{n})
+\frac{1}{n}\ln \nu(\{\Pi_{\Theta}(\boldsymbol {\theta}_0)\}\})
\leq \frac{1}{n}\ln \int_{\Theta} \exp (n\ell(\boldsymbol{\theta };\bar{\boldsymbol{x}}_{n}))\nu (d\boldsymbol{\theta }),$$
so letting $n \rightarrow \infty$ and using the continuity of $\ell(\Pi_{\Theta}(\boldsymbol {\theta}_0);\cdot)$, the required lower bound \begin{equation}\label{lb1}
\liminf _{n \rightarrow \infty }\frac{1}{n}\ln \int_{\Theta} \exp (n\ell(\boldsymbol{\theta };\bar{\boldsymbol{x}}_{n}))\nu (d\boldsymbol{\theta }) \geq \ell(\Pi_{\Theta}(\boldsymbol {\theta}_0); \boldsymbol{\mu}_0)=\kappa _{\Theta}^{\ast} ({\boldsymbol{\mu}}_0),
\end{equation}
is obtained. Next suppose instead that $\nu(\{\Pi_{\Theta}(\boldsymbol {\theta}_0)\})=0$; then, for any $\epsilon >0$ the ball
$B\left (\Pi_{\Theta}(\boldsymbol {\theta}_0),\epsilon \right )$ of radius $\epsilon$ around $\Pi_{\Theta}(\boldsymbol {\theta}_0)$ has a positive $\nu$-probability  for any $\epsilon>0$, so it holds
$$\inf _{\boldsymbol{\theta } \in \Theta \cap B\left (\Pi_{\Theta}(\boldsymbol {\theta}_0),\epsilon \right )}\ell(\boldsymbol{\theta };\bar{\boldsymbol{x}}_{n})
+\frac{1}{n}\ln \nu \left (B\left (\Pi_{\Theta}(\boldsymbol {\theta}_0) ,\epsilon \right )\right ) 
\leq \frac{1}{n}\ln \int_{\Theta} \exp (n\ell(\boldsymbol{\theta };\bar{\boldsymbol{x}}_{n}))\nu (d\boldsymbol{\theta }).$$
Then letting $n \rightarrow \infty$ first, and then $\epsilon \downarrow 0$, one gets
\begin{equation}\label{mancava}
\sup_{\epsilon > 0}\liminf _{n \rightarrow \infty }\inf _{\boldsymbol{\theta } \in \Theta \cap B(\Pi_{\Theta}(\boldsymbol {\theta}_0) ,\epsilon)}\ell(\boldsymbol{\theta };\bar{\boldsymbol{x}}_{n}) \leq \liminf _{n \rightarrow \infty }\frac{1}{n}\ln \int_{\Theta} \exp (n\ell(\boldsymbol{\theta };\bar{\boldsymbol{x}}_{n}))\nu (d\boldsymbol{\theta }) .
\end{equation}
The left hand side in display \eqref{mancava} cannot be smaller than $\kappa _{\Theta}^{\ast} ({\boldsymbol{\mu}}_0)=\ell(\Pi_{\Theta}(\boldsymbol {\theta}_0); \boldsymbol {\mu}_0)$. Suppose instead that for some $\delta>0$
$$\sup_{\epsilon > 0}\liminf _{n \rightarrow \infty }\inf _{\boldsymbol{\theta } \in \Theta \cap B(\Pi_{\Theta}(\boldsymbol {\theta}_0) ,\epsilon)}
\ell(\boldsymbol{\theta };\bar{\boldsymbol{x}}_{n})<\ell(\Pi_{\Theta}(\boldsymbol {\theta}_0);{\boldsymbol{\mu}}_{0})-\delta.$$
Then for any positive integer $m$ there exists $\boldsymbol{\theta}_m \in \Theta \cap B(\Pi_{\Theta}(\boldsymbol {\theta}_0) ,m^{-1})$ and an integer $n_m$ such that
\begin{equation}\label{assurdo}
\boldsymbol{\theta }_{m}\cdot \bar{\boldsymbol{x}}_{n_m} -\kappa(\boldsymbol{\theta }_{m})<
\Pi_{\Theta}(\boldsymbol {\theta}_0)\cdot {\boldsymbol{\mu}}_{0} -\kappa(\Pi_{\Theta}(\boldsymbol {\theta}_0))-\delta.
\end{equation}
Now $\boldsymbol{\theta }_{m}$ converges to $\Pi_{\Theta}(\boldsymbol {\theta}_0)$, $\kappa (\boldsymbol{\theta }_{m})$ converges to $\kappa (\Pi_{\Theta}(\boldsymbol {\theta}_0))$, and $n_m$ can be chosen to be increasing with $m$. As $m \to \infty$ we get the convergence of the left hand side of \eqref{assurdo} to $\Pi_{\Theta}(\boldsymbol {\theta}_0)\cdot {\boldsymbol{\mu}}_{0} -\kappa(\Pi_{\Theta}(\boldsymbol {\theta}_0))$, which is impossible since we assumed $\delta>0$.
So we have proved that \eqref{lb1} still holds, ending the proof.
\end{proof}

\paragraph{Proof of Theorem \ref{invsanofi}.} 
The proof of the upper bound consists in estimating the numerator of the Bayes' formula \eqref{bayes}. Choose a Borel set $B\subset \Theta$: then, with exactly the same argument of the previous lemma
\begin{equation}\label{suppartition}
\limsup_{n \to \infty} \frac{1}{n}\ln \int_{B} \exp (\ell(\boldsymbol{\theta };\bar{\boldsymbol{x}}_{n}))\nu (d\boldsymbol{\theta })\leq \lim_{n \to \infty}\kappa^{\ast}_{ B}(\bar {\boldsymbol{x}}_{n})=
\kappa^{\ast}_{ B}({\boldsymbol{\mu}}_0)=\sup_{\boldsymbol{\theta} \in B} \ell(\boldsymbol{\theta};{\boldsymbol{\mu}}_0)\leq \sup_{\boldsymbol{\theta} \in \bar {B}} \ell(\boldsymbol{\theta};{\boldsymbol{\mu}}_0).
\end{equation}
Indeed the supremum in \eqref{suppartition} is increased once the closure of $B$ in the relative topology of $\Theta$ is taken.

As far as the lower bound is concerned, let $O \cap \Theta $ be the interior of the measurable set $B$ in the relative topology 
of $\Theta$, $O$ being an open set in $\mathbb {R}^d$. Repeating the argument of the previous proof with any 
$\boldsymbol{\theta }_{\ast} \in O \cap \Theta$ replacing $\Pi_{\Theta}(\boldsymbol {\theta}_0)$, one arrives at
\begin{equation}\label{lb2}
\liminf _{n \rightarrow \infty }\frac{1}{n}\ln \int _{B\cap \Theta}\exp (n\ell(\boldsymbol{\theta };\bar{\boldsymbol{x}}_{n}))\nu (d\boldsymbol{\theta }) \geq \ell(\boldsymbol {\theta}_{\ast};\boldsymbol{\mu}_0).
\end{equation}
By combining Lemma \ref{th1}, \eqref{suppartition} and \eqref{lb2}, the LDP holds with rate function
$$
J(\boldsymbol {\theta})=\kappa _{\Theta}^{\ast} (\boldsymbol{\mu}_0)-\ell(\boldsymbol {\theta};\boldsymbol{\mu}_0)=
\big(\kappa^{\ast} (\boldsymbol{\mu}_0)-\ell(\boldsymbol {\theta};\boldsymbol{\mu}_0)\big)-\Big(\kappa^{\ast} (\boldsymbol{\mu}_0)-\kappa _{\Theta}^{\ast} (\boldsymbol{\mu}_0))\Big)
$$
which, by using \eqref{remember} and \eqref{eq:aggiunta}, is readily translated into the promised expression \eqref{invdiv}.\\
\ \\

Having proved our main result, we turn our attention to the characterization of the $m$-projection of the "true" parameter $\boldsymbol {\theta}_0$ on various classes of closed sets (in the relative topology of $\Theta_0$). This will turn out to be useful also for the discussion in the next section, devoted to of maximum likelihood estimation. We start with an equality which can be immediately checked, but it is extremely useful.

\begin{lemma}\label{threepoints} For any $\boldsymbol {\theta}_i \in \Theta_0, i=1,2,3$ the following holds
$$
D(\boldsymbol {\theta}_1||\boldsymbol {\theta}_3)=D(\boldsymbol {\theta}_1||\boldsymbol {\theta}_2) + D(\boldsymbol {\theta}_2||\boldsymbol {\theta}_3)
-(\boldsymbol {\theta}_3-\boldsymbol {\theta}_2) \cdot (\nabla \kappa (\boldsymbol {\theta}_1)-\nabla \kappa (\boldsymbol {\theta}_2)).
$$
\end{lemma}
From Lemma \ref{threepoints} we obtain firstly a necessary condition for the $m$-projection of $\boldsymbol {\theta}_0$ on the closed set $\Theta$, that we rephrase from Theorem 5.12 in Brown (1986), with the exception that we find convenient to define an \emph{admissible direction} to the set $\Theta$ at a point $\boldsymbol {\theta} \in \Theta$, as a limit point of a sequence $||\boldsymbol {\theta}_n-\boldsymbol {\theta}||^{-1}(\boldsymbol {\theta}_n-\boldsymbol {\theta})$, where $\{\boldsymbol {\theta}_n\}\subset \Theta$ is any sequence converging to $\boldsymbol {\theta}$, as $n \to \infty$. The \emph{cone of admissible directions} $A_{\boldsymbol {\theta}}(\Theta)$ to $\Theta$ at $\boldsymbol {\theta} \in \Theta$  is the closed convex cone generated by all the admissible directions. In case $\Theta$ is a differentiable manifold it coincides with the tangent space to $\Theta$ at $\boldsymbol {\theta}$. This notion allows a precise definition of the \emph{outward normal cone} 
$\nabla _{\Theta}(\boldsymbol {\theta})$ to $\Theta$ at $\theta \in \Theta$ as the polar cone to $A_{\boldsymbol {\theta}}(\Theta)$, that is the set of $\boldsymbol {v} \in \mathbb {R}^d$ such that
$$
\boldsymbol {w} \cdot \boldsymbol {v}\leq 0,\,\,\, \mbox{for all}\ \boldsymbol {w} \in A_{\boldsymbol {\theta}}(\Theta)
$$
(the same object defined by Brown in a slightly different way). When $\Theta$ is a differentiable manifold the outward normal cone $\nabla _{\Theta}(\boldsymbol {\theta})$ is the orthogonal complement of the tangent space to $\Theta$ at $\boldsymbol \theta$; when $\Theta$ is a convex set $\nabla _{\Theta}(\boldsymbol {\theta})$ is the normal cone to $\Theta$ at $\boldsymbol {\theta}$. The benefit of our definition is that, since the polarity operation on closed convex cones is an involution, the outward normal cone to  $\nabla _{\Theta}(\boldsymbol {\theta})$ at $\boldsymbol {\theta}$ is $A_{\boldsymbol {\theta}}(\Theta)$.

Now we give a necessary condition for the identification of $\Pi_{\Theta}(\boldsymbol {\theta}_0)$ for general $\Theta$.

\begin{proposition}\label{general}
Let $\Pi_{\Theta}(\boldsymbol {\theta}_0)\in \Theta$ be a $m$-projection of $\boldsymbol {\theta}_0 \in \Theta_0$ on the relatively closed subset $\Theta$
of $\Theta_0$. Then
\begin{equation}\label{first}
\nabla \kappa (\boldsymbol {\theta}_0)- \nabla \kappa (\Pi_{\Theta}(\boldsymbol {\theta}_0)) \in \nabla _{\Theta}(\Pi_{\Theta}(\boldsymbol {\theta}_0)).
\end{equation}
\end{proposition}
\begin{proof}
In Lemma \ref{threepoints} choose $\boldsymbol {\theta}_1=\boldsymbol {\theta}_0$, $\boldsymbol {\theta}_2=\Pi_{\Theta}(\boldsymbol {\theta}_0)$ and $\boldsymbol {\theta}_3=\boldsymbol {\theta}_n \in \Theta$ converging to $\Pi_{\Theta}(\boldsymbol {\theta}_0)$. Then, from \eqref{exponentialKL} the function $F(\boldsymbol {\theta})=D(\Pi_{\Theta}(\boldsymbol {\theta}_0)||\boldsymbol {\theta})$ has $F^{\prime}(\Pi_{\Theta}(\boldsymbol {\theta}_0))=0$, hence the result follows from the second order Taylor expansion of $F$ around $\Pi_{\Theta}(\boldsymbol {\theta}_0)$. 
\end{proof}

As the reader can imagine, when the set $\Theta$ is convex, since the function $D(\boldsymbol {\theta}_0||\boldsymbol {\theta})$ is strictly convex in $\boldsymbol {\theta} \in \Theta_0$, the $m$-projection is unique and it is possible to strengthen the previous proposition in the following way.

\begin{proposition}\label{projection} 
Let $\Theta$ be a relatively closed convex subset of $\Theta_0$, and $\boldsymbol {\theta}_0 \in \Theta_0$. Then $\Pi_{\Theta}(\boldsymbol {\theta}_0)\in \Theta$ is the $m$-projection of $\boldsymbol {\theta}_0$
on $\Theta$ if and only if one of the two following equivalent condition holds:
\begin{itemize}
\item
a. For any $\boldsymbol {\theta} \in \Theta$
\begin{equation}\label{serve}
(\boldsymbol {\theta}-\Pi_{\Theta}(\boldsymbol {\theta}_0)) \cdot (\nabla \kappa (\boldsymbol {\theta}_0)-\nabla \kappa (\Pi_{\Theta}(\boldsymbol {\theta}_0))) \leq 0;
\end{equation}
\item
b. For any $\boldsymbol {\theta} \in \Theta$
\begin{equation}\label{dispit}
D(\boldsymbol {\theta}_0||\boldsymbol {\theta})\geq 
D(\boldsymbol {\theta}_0||\Pi_{\Theta}(\boldsymbol {\theta}_0)) + 
D(\Pi_{\Theta}(\boldsymbol {\theta}_0)||\boldsymbol {\theta}).
\end{equation}
\end{itemize} 
If $\Theta$ is a \emph{linear subfamily}, i.e. the intersection with $\Theta_0$
of a linear manifold $L=\boldsymbol {a}+H$, $H$ being a proper subspace of $\mathbb {R}^d$,
then $\Pi_{\Theta}(\boldsymbol {\theta}_0)$ is the $m$-projection of $\boldsymbol {\theta}_0$
on $\Theta$ if and only if one of the two following equivalent condition holds:
\begin{itemize}
\item
c.
$$
\boldsymbol{h}\cdot(\nabla \kappa (\boldsymbol {\theta}_0)-\nabla \kappa (\Pi_{\Theta}(\boldsymbol {\theta}_0)))=0\ \mbox{for all}\ \boldsymbol{h}\in H,
$$
i.e. $\nabla \kappa (\boldsymbol {\theta}_0)-\nabla \kappa (\Pi_{\Theta}(\boldsymbol {\theta}_0))$ is orthogonal to $H$;
\item
d. the following Pythagorean identity holds
$$
D(\boldsymbol {\theta}_0||\boldsymbol {\theta})=
D(\boldsymbol {\theta}_0||\Pi_{\Theta}(\boldsymbol {\theta}_0)) + D(\Pi_{\Theta}(\boldsymbol {\theta}_0)||\boldsymbol {\theta})
\,\,\, \mbox{for all}\ \boldsymbol {\theta} \in \Theta.
$$
\end{itemize}
\end{proposition}

\begin{proof} Writing again Lemma \ref{threepoints} with $\boldsymbol {\theta}_1=\boldsymbol {\theta}_0$, $\boldsymbol {\theta}_2=\Pi_{\Theta}(\boldsymbol {\theta}_0)$ and $\boldsymbol {\theta}_3=\boldsymbol {\theta} \in \Theta$, we immediately see that the two conditions above are equivalent. 
It is equally straightforward to check that b. implies that $\Pi_{\Theta}(\boldsymbol {\theta}_0)$ is the $m$-projection
because it implies, for any $\boldsymbol {\theta} \in \Theta$,
$$
D(\boldsymbol {\theta}_0||\boldsymbol {\theta}) \geq D(\boldsymbol {\theta}_0||\Pi_{\Theta}(\boldsymbol {\theta}_0)).
$$
On the other hand, thanks to Lemma \ref{threepoints}, condition a. is nothing but \eqref{first} for $\Theta$ convex. 
The equivalent conditions c. and d. are immediately obtained from a. and b. for the case of a linear subfamily.
\end{proof}

Essentially with the same proof, Proposition \ref{projection} (and the forthcoming corollary) can be extended to the case of a \emph{quasi-convex set} $\Theta$, that is a set with the property that, whenever $\boldsymbol {\theta}, \boldsymbol {\theta}_0 \in \Theta$, then $\boldsymbol {\theta}-\boldsymbol {\theta}_0$  is an admissible direction to $\boldsymbol {\theta}_0$ at $\Theta$ (see Miao and Hahn, 1997).

The first mention of the Pythagorean identity in item d. of Proposition \ref{projection} seems to be Chentsov (1966); other relevant references are Simon (1973) and Csiszar (1975). Incidentally, an alternative proof of the identity can be given exploiting the fact that a linear subfamily $L=a+H$ of a full natural exponential family is by itself a general exponential family with a sufficient statistics given by the projection onto $H$  (see Barndorff-Nielsen, 1978, Section 8.2). 

Taking Proposition \ref{projection} into account, Theorem \ref{invsanofi} yields the following result in the case of a linear subfamily or, more generally, a convex subfamily. In the first case, even if the model is misspecified and the true value $\boldsymbol {\theta}_0 \notin \Theta$, the asymptotic behaviour of the sequence of posterior distribution is equivalent, as far as the Large Deviation rate is concerned, to a well-specified family with the true value $\Pi_{\Theta}(\boldsymbol {\theta}_0) \in \Theta$. At the contrary, for a convex but misspecified family, the rate function is generally larger compared to the well-specified regime: the misspecification is felt at the level of large deviation, as it will be shown by the next example.

\begin{corollary}\label{invsanofi2} 
If the sequence  $\{\bar{\boldsymbol{x}}_{n}, n\in \mathbb{N}\}$ converges to some $\boldsymbol{\mu}_0 \in M_0$, with $\boldsymbol {\mu}_0=\nabla \kappa (\boldsymbol {\theta}_0)$, the sequence of probability measures $\{\pi_{n}(\cdot|\bar {\boldsymbol{x}}_{n}), n \in \mathbb {N}\}$ satisfies a LDP
\begin {itemize}
\item a.
with rate function 
$$J(\boldsymbol {\theta})=D(\Pi_{\Theta}(\boldsymbol {\theta}_0)||\boldsymbol {\theta}),$$ if $\nu$ is supported by a linear subfamily $\Theta$;
\item
b. with rate function 
\begin{equation}\label{pure}
J(\boldsymbol {\theta})=D(\Pi_{\Theta}(\boldsymbol {\theta}_0)||\boldsymbol {\theta})-(\boldsymbol {\theta}-\Pi_{\Theta}(\boldsymbol {\theta}_0)) \cdot (\nabla \kappa (\boldsymbol {\theta}_0)-\nabla \kappa (\Pi_{\Theta}(\boldsymbol {\theta}_0))
\geq D(\Pi_{\Theta}(\boldsymbol {\theta}_0)||\boldsymbol {\theta}),
\end{equation}
if $\nu$ is supported by a convex subfamily $\Theta$.
\end{itemize}
\end{corollary}
\begin{proof}
The result is a straightforward consequence of Proposition \ref{projection}. Notice that the inequality sign \eqref{pure} follows from that in \eqref{serve}.
\end{proof}

\begin{example}\rm
The Hardy-Weinberg family of distributions, in its simplest form with two alleles (see e.g. Barndorff-Nielsen (1978), Example 8.10), is a subfamily 
of the saturated family of all distributions over $3$ outcomes (the genotypes $AA$, $aa$ and $Aa$), coded with the $3$ vectors in the plane $\mathbf {e}_1=(1,0), \mathbf {e}_2=(0,1)$ and $\mathbf {0}=(0,0)$, respectively. Indeed the latter family is represented as the natural exponential family $\{P_{\boldsymbol {\theta}}, \boldsymbol {\theta} \in \Theta_0\}$ generated by 
$$\lambda=\frac {1}{2}\delta_{\mathbf {0}}+\frac {1}{4}\delta_{\mathbf {e}_1}+\frac {1}{4}\delta_{\mathbf {e}_2},$$
whose cumulant generating function is
$$\kappa (\boldsymbol {\theta})=\log \int \exp \{\boldsymbol {\theta} \cdot \mathbf {x}\}\lambda (d\mathbf {x})=\log (2+e^{\theta_1}+e^{\theta_2})-2\log {2},\,\,\,\,\, \boldsymbol {\theta}=(\theta_1, \theta_2) \in \mathbb {R}^2.$$
The probabilities that the element of the exponential family $P_{\boldsymbol {\theta}}$ defined in \eqref{family} assigns to the three outcomes are
\begin{equation}\label{outcomes}
P_{\boldsymbol{\theta }}(\mathbf {e}_i)=\frac {\partial \kappa}{\partial \theta_i}(\theta_1, \theta_2)=\frac {e^{\theta_i}}{2+e^{\theta_1}+e^{\theta_2}} 
\,\,\,\ \mbox{for}\ i=1,2, \,\,\,\ \mbox{and}\ \,\,\,\ P_{\boldsymbol{\theta }}(\mathbf{0})=\frac {2}{2+e^{\theta_1}+e^{\theta_2}}.
\end{equation}
It is clear that $C(\lambda)$ is the triangle with vertices $\mathbf {0}=(0,0), \mathbf {e}_1=(1,0), \mathbf {e}_2=(0,1)$. The relation \eqref{outcomes} can be easily inverted: if $(x,y) \in M_0=int(C(\lambda))$, meaning that $x>0, y>0, 1-x-y>0$, then, setting the probabilities of the outcomes in \eqref{outcomes} equal to $x$, $y$ and $1-x-y$, one can check that
\begin{equation}\label{nathw}
\theta_1=\ln 2+\ln x -\ln (1-x-y)=\frac {\partial \kappa^{\ast}}{\partial x}(x,y),\,\,\, 
\theta_2=\ln 2+\ln y -\ln (1-x-y)=\frac {\partial \kappa^{\ast}}{\partial y}(x,y),
\end{equation}
where $\kappa^{\ast}$ is defined in \eqref{supachie}.

The Hardy-Weinberg subfamily assumes that the probabilities in \eqref{outcomes} arise from a binomial distribution with $2$ trials, where 
$\mathbf {e}_1$ corresponds to two successes, $\mathbf {e}_2$ corresponds to two failures, and
$\mathbf{0}$ corresponds to one success and one failure. This identifies the subset
$$M=\{(x,y)\in M_0:xy=(1-x-y)^2\},$$
which in term of the natural parameters becomes
$$
\Theta=\nabla \kappa^{\ast}(M)=\{(\theta_1, \theta_2):\theta_1+\theta_2=0\},
$$
assumed to be the support of the prior distribution $\nu$. For $\boldsymbol {\theta}=(\theta_1, -\theta_1) \in \Theta$, setting
$$
p_A(\theta_1)=\frac{e^{\theta_1/2}}{\sqrt{2(1+\cosh(\theta_1))}}\,\,\,\,\mbox{and}\,\,\,\, p_a(\theta_1)=1-p_A(\theta_1)=\frac{e^{-\theta_1/2}}{\sqrt{2(1+\cosh(\theta_1))}},
$$
the probabilities of the three outcomes become
$$P_{\boldsymbol{\theta }}(\mathbf {e}_1)=p_A^2(\theta_1),\,\,\,\,
P_{\boldsymbol{\theta }}(\mathbf {e}_2)=p_a^2(\theta_1)\,\,\,\,\mbox{and}\,\,\,\,
P_{\boldsymbol{\theta }}(\mathbf {0})=2p_A(\theta_1)p_a(\theta_1).$$

For $\boldsymbol {\mu}_0=(x_0,y_0)\in M_0$, let $P_{\boldsymbol {\theta}_0}$ be the element of the saturated model such that $P_{\boldsymbol {\theta}_0}(\mathbf {e}_1)=x_0$ and $P_{\boldsymbol {\theta}_0}(\mathbf {e}_2)=y_0$, obtained through \eqref{nathw}.
By applying Proposition \ref{projection} the $m$-projection $\Pi_{\Theta}(\boldsymbol {\theta}_0)$ of $\boldsymbol {\theta}_0$ on $\Theta$ 
is characterized by the fact that $(\hat x, \hat y)=\nabla\kappa(\Pi_{\Theta}(\boldsymbol {\theta}_0))$ is obtained by joining $\boldsymbol {\mu}_0$ with $M$ through a $45^{\circ}$ line. It is easily obtained that 
\begin{equation}\label{hat}
\hat x =\left(\frac {1+x_0-y_0}{2}\right)^2, \hat y =\left(\frac {1+y_0-x_0}{2}\right)^2, 1-\hat x-\hat y=\frac {(1+x_0-y_0)(1+y_0-x_0)}{2}=\frac {1-(x_0-y_0)^2}{2}.
\end{equation}
This has a biological meaning: the frequencies $(\hat x, \hat y, 1-\hat x-\hat y)$ express the so-called Hardy-Weinberg equilibrium, achieved for the first (infinitely large) generation of alleles starting in the generation of parents with the frequencies $(x_0, y_0, 1-x_0-y_0)$ for $AA$, $aa$ and $Aa$, respectively.

By Theorem \ref{invsanofi}, if the sequence  $\{\bar{\boldsymbol{x}}_{n}, n\in \mathbb{N}\}$ converges to $\boldsymbol {\mu}_0=(x_0,y_0)$ as $n\to\infty$, the sequence of probability measures $\{\pi_{n}(\cdot|\bar {\boldsymbol{x}}_{n}), n \in \mathbb {N}\}$ on $\Theta$ satisfies a LDP with rate function 
\begin{multline*}
J_1(\theta_1 ):= J(\theta_1, -\theta_1)=\hat x \log \frac{\hat x}{p_A(\theta_1)^2}+ \hat y \log \frac{\hat y}{p_a(\theta_1)^2}+
(1-\hat x - \hat y) \log \frac {1-\hat x - \hat y}{2p_A(\theta_1)p_a(\theta_1)}\\
=(1+x_0-y_0) \log \frac {1+x_0-y_0}{2p_A(\theta_1)}+(1+y_0-x_0) \log \frac {1+y_0-x_0}{2p_a(\theta_1)}.
\end{multline*}
Not surprisingly, this is twice the KL divergence of a Bernoulli distribution, with parameter $\frac {1+x_0-y_0}{2}$, with respect to a Bernoulli distribution, with parameter $p_A(\theta_1)$. This is unchanged if $x_0$ and $y_0$ are replaced by $\hat x$ and $\hat y$, which illustrates case a. in Corollary \ref{invsanofi2}.

In order to explore case b., let us set the model to be unfavourable for the 
heterozygous allele $Aa$ compared to the Hardy-Weinberg model, meaning that the mean value is assumed to belong to the subset  
$$
M^+=\{(x,y) \in M_0: xy\geq(1-x-y)^2 \}
$$
whose image in the natural parameter space is 
$$
\Theta^+=\nabla \kappa^{\ast}(M^+)=\{(\theta_1, \theta_2):\theta_1+\theta_2\geq 0\}.
$$
In this case, if $\boldsymbol {\mu}_0=(x_0,y_0) \in M^+$, the $m$-projection of $\boldsymbol {\theta}_0$ on $\Theta^+$ is $\boldsymbol {\theta}_0$ itself, so the rate function $J(\boldsymbol \theta)$ is the KL divergence $D(\boldsymbol {\theta}_0||\boldsymbol {\theta})$, for any $\boldsymbol {\theta}\in \Theta^+$. But if $\boldsymbol {\mu}_0=(x_0,y_0)\in M_0\setminus M^+$ the $m$-projection of $\boldsymbol {\theta}_0$ on $\Theta^+$ coincides with that on $\Theta$, the parameter space of the Hardy-Weinberg family, with the same mean vector $(\hat x, \hat y)$ written in \eqref{hat}. The rate function can be obtained from the expression \eqref{pure}.

The previous example actually contains an extension of the results of Ganesh and O'Connell (1999) for the multivariate Bernoulli framework, in that  misspecification is allowed.
\end{example}

\section{Large deviation principles for the MLE}\label{sec:LDP-MLE}

This section is devoted to the discussion of the frequentist counterpart of Theorem \ref{invsanofi}, namely the LDP for MLE in a curved exponential family. Our goal is to emphasize the link with information geometry and the relation with the "parametric" inverse Sanov theorem, established by 
Theorem \ref{invsanofi} for well-specified models. For this reason we are going to made rather heavy simplifying assumptions. At the same time the assumption that the natural parameter space of the exponential family \eqref{family} is open can be relaxed, which will turn out to be useful in the next section. It is enough to assume that $\kappa$ is \emph{steep}, which means that $\Theta_0$ has a non-empty interior (over which the formula \eqref{grad} for the gradient of $\kappa$ is always guaranteed), and, for any $\boldsymbol {\theta} \in \partial {\Theta_0}$ and any sequence $\{ \boldsymbol {\theta}_n, n \in \mathbb{N} \} \in int (\Theta)$ converging to $\boldsymbol {\theta}$, the sequence $\{ ||\nabla (\boldsymbol {\theta}_n)||, n \in \mathbb{N} \}$ diverges, as $n \to \infty$. In the steep case, the fact that  $\nabla \kappa$ and $\nabla \kappa^{\ast}$ are inverse to each other, between $M_0$ and $int(\Theta_0)$, still holds true (Barndorff-Nielsen, Theorems 5.33 and 9.2).

Now let \textbf{$\boldsymbol{x}_{1} ,\boldsymbol{x}_{2} ,\boldsymbol{\ldots } ,\boldsymbol{x}_{\boldsymbol{n}}$} be an i.i.d. sample drawn from $P_{\boldsymbol{\theta }_0}$ belonging to the family \eqref{family}, where $\boldsymbol{\theta }_0 \in int(\Theta_0)$. Our next assumption is that, under $P_{\boldsymbol{\theta }_0}$, the sufficient statistics $\bar{\boldsymbol{x}}_{n}$ takes values in $M_0$ with probability $1$ (this condition clearly holds irrespectively of $\boldsymbol{\theta }_0$). As long as $\bar{\boldsymbol{x}}_{n}$ takes values in $M_0$ the log-likelihood function
$\ell(\boldsymbol {\theta};\bar{\boldsymbol{x}}_{n})$  has an unconstrained maximizer $\nabla \kappa^{\ast}(\bar{\boldsymbol{x}}_{n})$ in $\boldsymbol {\theta} \in \Theta_0$. We have already seen that the maximization of the log-likelihood within a subfamily indexed by $\Theta \subset \Theta_0$ is equivalent to the minimization of $D(\nabla\kappa^{\ast}(\bar{\boldsymbol{x}}_{n})||\boldsymbol {\theta})$ with respect to $\boldsymbol {\theta} \in \Theta$,  i.e. maximum likelihood estimates are $m$-projections of $\nabla \kappa^{\ast}(\bar{\boldsymbol{x}}_{n})$ on the subfamily $\Theta$, which exist as long as $\Theta$ is relatively closed in $\Theta_0$. By a Maximum Likelihood Estimator \emph{constrained} to $\Theta$ we mean a Borel function $\varphi_{\Theta}$ of $\bar{\boldsymbol{x}}_{n}\in M_0$ taking values in the set of $m$-projections of $\nabla \kappa^{\ast}(\bar{\boldsymbol{x}}_{n})$ on $\Theta$.

In order to state the next result, let us introduce the set
$$
T_{\boldsymbol {\theta}}(\Theta)=\{ \boldsymbol {\theta}'\in int(\Theta_0): \varphi_{\Theta} (\nabla \kappa (\boldsymbol {\theta}'))=\boldsymbol {\theta}\}.
$$ 
Then the following result holds.

\begin{theorem}\label{mle} Suppose that $\boldsymbol {x}_1,\ldots, \boldsymbol {x}_n$ is an i.i.d. sample drawn from
$P_{\boldsymbol{\theta }_0}$ and the sample mean $\bar{\boldsymbol{x}}_{n}$ takes values in $M_0$ almost surely. Let $\Theta \subset \Theta_0$ be a relatively closed subset of $\Theta_0$ and suppose that there exists a continuous function $\varphi_{\Theta}: M_0 \to \Theta \cap int(\Theta_0)$ which is a MLE constrained to $\Theta$. Then the sequence of laws of $\{\varphi_{\Theta}(\bar{\boldsymbol{x}}_{n}), n \in \mathbb{N}\}$ satisfies a LDP with a rate function 
\begin{equation}\label{cracontr}
I(\boldsymbol{\theta})=
\inf \{D(\boldsymbol{\theta}'||\boldsymbol{\theta }_0): \boldsymbol{\theta' } \in T_{\boldsymbol {\theta}}(\Theta)\},\,\,\, \boldsymbol{\theta} \in \Theta \cap int(\Theta_0).
\end{equation}
\end{theorem}

\begin{proof} 
By Cram\'er's theorem (see e.g. Theorem 2.2.30 in Dembo and Zeitouni (1998)), the sample mean of i.i.d. random variables
drawn from $P_{\boldsymbol{\theta }_0}$ satisfies a LDP with rate function
$$\iota (\boldsymbol {t})=\sup_{\boldsymbol{\theta}\in \mathbb {R}^d}\left\{\boldsymbol{t}\cdot\boldsymbol{\theta}
-\log\int e^{(\boldsymbol {\theta}_0+\boldsymbol{\theta})\cdot \boldsymbol {x}-\kappa (\boldsymbol {\theta}_0)}\lambda(d\boldsymbol {x})\right\}
=\sup_{\boldsymbol{\theta}\in \mathbb {R}^d}\left\{\boldsymbol{t}\cdot\boldsymbol{\theta}
-\kappa (\boldsymbol {\theta}+\boldsymbol {\theta}_0)+\kappa (\boldsymbol {\theta}_0)\right\}.$$
Therefore, by (\ref{achieve}) and a straightforward change of variable, for $\boldsymbol {t} \in M_0$
we have
$$\iota (\boldsymbol {t})=\kappa^{\ast} (\boldsymbol {t})-l(\boldsymbol {\theta}_0;\boldsymbol {t})
=D(\nabla \kappa^{\ast} (\boldsymbol{t})||\boldsymbol{\theta }_0).$$
Finally, since $\varphi_{\Theta}$ is continuous, we use the contraction principle (see e.g. Dembo and Zeitouni (1998), Theorem 4.2.1) to establish the
LDP for the sequence of laws of $\{\varphi_{\Theta}(\bar{\boldsymbol{x}}_{n}), n \in \mathbb{N}\}$ with rate function $I$ 
defined by
$$
I(\boldsymbol {\theta})=
\inf \{\iota (\boldsymbol {t}): \boldsymbol {t} \in M_0, \varphi_{\Theta}(\boldsymbol {t})=\boldsymbol{\theta}\}
=\inf \{ \iota (\nabla \kappa (\boldsymbol {\theta}')): \boldsymbol {\theta}' \in 
T_{\boldsymbol {\theta}}(\Theta)\},\ \boldsymbol{\theta} \in \Theta \cap int(\Theta_0)$$
which is easily seen to coincide with \eqref{cracontr} since $\nabla \kappa$ and $\nabla \kappa^{\ast}$ are inverse to each other between $M_0$ and $int(\Theta_0)$.
\end{proof}

As a consequence of the previous result, under the conditions stated therein, the rate function $I(\boldsymbol {\theta})$, for any $\boldsymbol {\theta} \in \Theta$, can be computed by finding a natural parameter value $\boldsymbol {\theta}^{\ast}\in \Theta$ which minimizes $D(\boldsymbol {\theta}'||\boldsymbol {\theta}_0)$ with $\boldsymbol {\theta}'$ belonging to $T_{\boldsymbol {\theta}}(\Theta)$, that is finding an $e$-projection ${\boldsymbol \theta^{\ast}}=\Pi^{T_{\boldsymbol {\theta}}(\Theta)}(\boldsymbol{\theta}_0)$ of ${\boldsymbol {\theta}_0}$ on $T_{\boldsymbol {\theta}}(\Theta)$. In the well-specified case, i.e. $\boldsymbol {\theta}_0 \in \Theta$, if for any $\boldsymbol {\theta} \in \Theta$ it holds
\begin{equation}\label{pita}
D(\boldsymbol {\theta}||\boldsymbol {\theta}_0)\leq D(\boldsymbol {\theta}'||\boldsymbol {\theta}_0),\,\,\, \mbox{for all}\ 
\boldsymbol {\theta}' \in T_{\boldsymbol {\theta}}(\Theta),
\end{equation}
then $\boldsymbol{\theta}=
\Pi^{T_{\boldsymbol {\theta}}(\Theta)}(\boldsymbol{\theta}_0)$ and
\emph{the "parametric" Sanov theorem holds}, i.e. $I(\boldsymbol {\theta})=D(\boldsymbol{\theta } ||\boldsymbol{\theta }_0)$. 

Note that, when $\Theta=\Theta_0$ and $\bar {\boldsymbol{x}}_{n}$ takes values in $M_0$ with probability 1, Theorem \ref{mle} can be applied and 
the property \eqref{pita} is trivially verified, since the set $T_{\boldsymbol {\theta}}(\Theta)$ reduces to the point $\{\boldsymbol {\theta}\}$. The following result extends the class of families for which \eqref{pita} is true to convex families; a closely related result result can be found in Theorem 2.2 in Kester and Kallenberg (1986).

In order to use in the following discussion the results in Section \ref{sec:preliminary}, we return for the rest of this section to the assumption that $\Theta_0$ is open.
 
\begin{theorem}\label{pyta1} 
Let $\Theta$ be a convex relatively closed subset of the open set $\Theta_0$, and $\boldsymbol {\theta}_0 \in \Theta$.  Suppose that under i.i.d. sampling from $P_{\boldsymbol {\theta}_0}$, the sample mean $\bar{\boldsymbol{x}}_{n}$ takes values in $M_0$ almost surely. Then there is a unique MLE $\varphi_{\Theta}$ defined in $M_0$ and the sequence of laws of $\{\varphi_{\Theta}(\bar{\boldsymbol{x}}_{n}), n \in \mathbb{N}\}$ satisfies the LDP with rate function
$$I(\boldsymbol {\theta})=D(\boldsymbol {\theta}||\boldsymbol {\theta}_0),\,\,\,\, \boldsymbol \theta \in \Theta.$$
\end{theorem}
\begin{proof}
Let us prove that the assumptions of Theorem \ref{mle} are verified. We have already seen that a MLE $\varphi_{\Theta}(\boldsymbol {t})$, for $\boldsymbol t \in M_0$, coincides a $m$-projection of $\boldsymbol {\theta}_0=\nabla \kappa^{\ast}(\boldsymbol t)$ on $\Theta$. When $\Theta$ is convex and relatively closed in $\Theta_0$, the function $\varphi_{\Theta}$ is infinitely differentiable, since it is the smooth gradient of the convex lower semi-continuous function $\kappa^{\ast}_{\Theta}$ defined in \eqref{maximum} (see Barndorff-Nielsen, 1978, Theorem 9.18). Thus the continuity assumption of Theorem \ref{mle} is guaranteed.

It remains to prove the expression for the rate function.
For the reader's convenience we rewrite \eqref{dispit} as
$$
D(\boldsymbol {\theta}^{\prime}||\boldsymbol {\theta}_0)\geq 
D(\boldsymbol {\theta}^{\prime}||\boldsymbol {\theta}) + D(\boldsymbol {\theta}||\boldsymbol {\theta}_0)
$$
where $\boldsymbol {\theta}_0,\boldsymbol {\theta} \in \Theta$ and $\boldsymbol {\theta}^{\prime}\in T_{\boldsymbol {\theta}}(\Theta)$. We immediately see that $\boldsymbol {\theta}=\Pi^{T_{\boldsymbol {\theta}}(\Theta)}(\boldsymbol{\theta}_0)$, from which the expression given for the rate function follows.
\end{proof}

Finally, in order to investigate the possibility of extending the previous result further, we make more explicit the fact that the dual structure of $m$ and $e$-projections comes is a result of the conjugacy between $\kappa$ and $\kappa^{\ast}$. Having this in mind, for $\boldsymbol{\theta} \in \mathbb {R}^d$, we define the function
$$
\ell^{\ast}(\boldsymbol{t}; \boldsymbol{\theta})=\boldsymbol{t}\cdot \boldsymbol{\theta}-\kappa^{\ast} (\boldsymbol{t}),\ \boldsymbol{t} \in \mathbb {R}^d.
$$ 
Again, when $\boldsymbol {t}$ and $\boldsymbol {\theta}$ are a conjugate pair, which is equivalent to $\boldsymbol {\theta} \in \Theta_0$ and $\boldsymbol {t}= \nabla \kappa (\boldsymbol {\theta})\in M_0$,
\begin{equation}\label{pressure}
\kappa (\boldsymbol{\theta})=\nabla \kappa (\boldsymbol {\theta})\cdot \boldsymbol{\theta}-\kappa^{\ast} (\nabla \kappa (\boldsymbol {\theta}))
=\ell^{\ast}(\nabla \kappa (\boldsymbol {\theta}); \boldsymbol {\theta})=\sup_{t \in \mathbb {R}^d}\ell^{\ast}(\boldsymbol{t}; \boldsymbol{\theta}).
\end{equation}

This relation allows to express the KL divergence between elements of the exponential family \eqref{family} in terms of the functions $\kappa$ and $\kappa^{\ast}$. We have from \eqref{exponentialKL}
\begin{equation}\label{againdiv}
D(\boldsymbol{\theta}_0||\boldsymbol{\theta})=\ell(\boldsymbol{\theta}_0;\nabla \kappa (\boldsymbol {\theta}_0))-\ell(\boldsymbol{\theta};\nabla \kappa (\boldsymbol {\theta}_0))=\kappa^{\ast}(\nabla \kappa (\boldsymbol {\theta}_0))-\ell(\boldsymbol{\theta};\nabla \kappa (\boldsymbol {\theta}_0)),
\end{equation}
but also, from \eqref{pressure},
\begin{equation}\label{againdiv2}
D(\boldsymbol{\theta}||\boldsymbol{\theta}_0)
=\kappa^{\ast}(\nabla \kappa (\boldsymbol {\theta}))-\ell(\boldsymbol{\theta}_0;\nabla \kappa (\boldsymbol {\theta}))
=\kappa (\boldsymbol {\theta}_0)-\ell^{\ast}(\nabla \kappa (\boldsymbol {\theta}); \boldsymbol {\theta}_0),
\end{equation}
for $\boldsymbol {\theta}, \boldsymbol {\theta}_0 \in \Theta_0$. As a consequence, any maximizer of $\ell^{\ast}(\nabla \kappa (\boldsymbol {\theta}); \boldsymbol {\theta}_0)$ for $\boldsymbol {\theta}\in\Theta$ is a $e$-projection
of $\boldsymbol {\theta}_0$ on $\Theta$.

\begin{theorem}\label{eproj}
Let $\Theta \subset \Theta_0$ be a relatively closed subset of the open set $\Theta_0$. Define
$$
B_{\boldsymbol {\theta}}:=(\nabla \kappa (\boldsymbol {\theta})+ \nabla_{\boldsymbol {\theta}}(\Theta))\cap M_0,\ \boldsymbol {\theta}\in\Theta,
$$
where $\nabla_{\boldsymbol {\theta}}(\Theta)$ is the outward normal cone to $\Theta$ at $\boldsymbol {\theta}$. Assume that $\{B_{\boldsymbol {\theta}},\boldsymbol {\theta}\in\Theta\}$ is a partition of $M_0$ and the function $\varphi:M_0\to\Theta$ defined by
$$\varphi(\boldsymbol {\mu})=\boldsymbol {\theta}
\ \Longleftrightarrow \ \boldsymbol {\mu}\in B_{\boldsymbol {\theta}},$$
is continuous. Moreover assume that, for some $\boldsymbol{\theta }_0\in\Theta$, whenever $\boldsymbol {x}_1,\ldots, \boldsymbol {x}_n$ is an i.i.d. sample drawn from
$P_{\boldsymbol{\theta }_0}$, the sample mean $\bar{\boldsymbol{x}}_{n}$ takes values in $M_0$ almost surely. Then the following statements hold.
\begin{itemize}
\item a. The function $\varphi_{\Theta}=\varphi$ satisfies the assumptions in Theorem \ref{mle}, with $T_{\boldsymbol {\theta}}(\Theta)=\nabla\kappa^*(B_{\boldsymbol {\theta}})$.
\item b. The inequality
$$
D(\boldsymbol {\theta}||\boldsymbol {\theta}_0)\leq D(\boldsymbol {\theta}'||\boldsymbol {\theta}_0),\,\,\, \mbox{for all}\ 
\boldsymbol {\theta}' \in T_{\boldsymbol {\theta}}(\Theta),
$$
holds if and only if $\boldsymbol {\theta}_0 \in (\boldsymbol {\theta}+
A_{\boldsymbol {\theta}}(\Theta))\cap\Theta$, where $A_{\boldsymbol {\theta}}(\Theta)$ is the cone of admissible directions to $\boldsymbol {\theta}$ in $\Theta$.
\item c. The sequence of laws of $\{\varphi_{\Theta}(\bar{\boldsymbol{x}}_{n}), n \in \mathbb{N}\}$ satisfies the LDP with rate function
\begin{equation}\label{punti}
I (\boldsymbol \theta)=D(\boldsymbol {\theta}||\boldsymbol {\theta}_0)
\end{equation}
if and only if $\boldsymbol {\theta}_0 \in (\boldsymbol {\theta}+
A_{\boldsymbol {\theta}}(\Theta))\cap\Theta$, otherwise it holds 
$I (\boldsymbol \theta)<D(\boldsymbol {\theta}||\boldsymbol {\theta}_0)$.
\end{itemize}
\end{theorem}
\begin{proof} For proving a. we have to check that $\boldsymbol {\theta}\in \Theta$ is a maximizer of the log-likelihood function $\ell(\cdot;\boldsymbol {\mu})$ constrained to $\Theta$ as long as $\boldsymbol {\mu} \in B_{\boldsymbol {\theta}}$. This is due to Proposition \ref{general}, which guarantees the necessity of this condition, and the assumption that $\{B_{\boldsymbol {\theta}},\boldsymbol {\theta}\in\Theta\}$ is a partition of $M_0$, which implies its sufficiency. The characterization of $T_{\boldsymbol {\theta}}(\Theta)$ is then immediately obtained. 

For proving b. we observe that, by exchanging $\kappa$ with $\kappa^{\ast}$ and
$\ell$ with $\ell^{\ast}$, a characterization analogous to \eqref{serve} is obtained for $e$-projections on sets whose $\nabla \kappa$ image is closed in $M_0$ and convex (see Theorem 6.13 in Brown, 1986). Since $B_{\boldsymbol {\theta}}$ is a closed convex cone and $\nabla \kappa^{\ast}$ a homeomorphism, this implies that the $e$-projection of any $\boldsymbol {\theta}_0\in \Theta$ on $T_{\boldsymbol {\theta}}(\Theta)$  is $\boldsymbol \theta$ itself if and only if $\boldsymbol {\theta}_0-\boldsymbol {\theta}$ lies in the polar cone to 
$T_{\boldsymbol {\theta}}(\Theta)$ that, by the involutive property of the polarity operation on closed convex cones, coincides with $A_{\boldsymbol {\theta}}(\Theta)$.

The statement c. follows from part b. and Theorem \ref{mle}.
\end{proof}

Here is an example in which the equality \eqref{punti} is never achieved, except for the trival case $\boldsymbol {\theta} = \boldsymbol {\theta}_0$.

\begin{example}\label{cvequalone}\rm
The two-parameter Gaussian exponential family is a natural exponential family with respect to the reference measure $\lambda$ which is the image of the Lebesgue measure on the real line under the mapping $x \mapsto (x, -x^2/2)$. Denoting by $\mu$ and $\sigma^2$ the mean and the variance of the density, the natural parameters $(\theta_1, \theta_2)$ are then
$$
\theta_1=\frac {\mu}{\sigma^2} \in \mathbb{R},\,\,\,\, \theta_2=\frac {1}{\sigma^2} >0,
$$ 
the cumulant generating function is
$$
\kappa (\theta_1, \theta_2)=\frac {1}{2}\left (\log {(2\pi)}- \log {\theta_2}+\frac {\theta_1^2}{\theta_2}\right ),
$$
and the gradient of $\kappa$ is given by
$$\frac {\partial \kappa}{\partial \theta_1}=\frac {\theta_1}{\theta_2},\,\,\, \frac {\partial \kappa}{\partial \theta_2}=-\frac {\theta_1^2}{2\theta_2^2}-\frac {1}{2\theta_2}.$$
Since $\lambda$ is supported by the graph of a concave function, the set $M_0=int (C(\lambda))$ is the subset of the plane below this graph. It is clear that, if $X_1,\ldots,X_n$ is a i.i.d. sample drawn from a Gaussian law and
$$\bar{X}_n=\frac{\sum_{i=1}^nX_i}{n},\quad \overline{X^2}_n=\frac{\sum_{i=1}^nX_i^2}{n},$$
then the vector of sufficient statistics $(\bar{X}_n,-\frac{1}{2}\overline{X^2}_n)\in M_0$ 
unless all the elements of the sample are equal, which clearly happens with probability $0$. 

The subfamily of laws with mean equal to the standard deviation corresponds to the following curve in the natural parameter space
$$
\Theta=\left \{(\theta_1,\theta_2): \theta_2=\theta_1^2,\,\,\,\, \theta_1>0\right \}.
$$
A point in this set can be written as $\boldsymbol {v}(\theta_1)=(\theta_1,\theta_1^2)$, with $\theta_1>0$ and
$$
\frac {\partial \kappa}{\partial \theta_1}(\boldsymbol {v}(\theta_1))=\frac {1}{\theta_1},\,\,\,\,\, \frac {\partial \kappa}{\partial \theta_2}(\boldsymbol {v}(\theta_1))=-\frac {1}{\theta_1^2},
$$
and setting $\boldsymbol {w}(\mu_1)=(\mu_1, -\mu_1^2)$, this is rewritten as
$$
\nabla \kappa (\boldsymbol {v}(\theta_1))=\boldsymbol {w} (\theta_1^{-1}), \theta_1>0,
$$
from which the image of $\Theta$ under the mapping $\nabla \kappa$ is the set
$$
M=\left \{(\mu_1,\mu_2): \mu_2=-\mu_1^2,\,\,\,\, \mu_1>0\right \}.
$$

According to Proposition \ref{general} a necessary condition for $\theta_1$ to be a maximum point of the log-likelihood function
$$
\ell(\boldsymbol {v}(\theta_1);\boldsymbol {\mu})=\boldsymbol {v}(\theta_1) \cdot \boldsymbol {\mu} - \kappa (\boldsymbol {v}(\theta_1)), \,\,\,\,\theta_1>0
$$
is that $\boldsymbol {\mu}-\boldsymbol {w} (\theta_1^{-1})$ is normal to the tangent vector $(1,2\theta_1)$ 
to $\Theta$ at $\boldsymbol {v}(\theta_1)$, that spans the cone of admissible directions $A_{\boldsymbol {v}(\theta_1)}(\Theta)$. Correspondingly $(-2\theta_1,1)$ spans the outward normal cone $\nabla_{\boldsymbol {v}(\theta_1)}(\Theta)$, so that
$$
B_{\boldsymbol {v}(\theta_1)}=\left\{\boldsymbol {\mu}\in M_0:\boldsymbol {\mu}=\boldsymbol {w} (\theta_1^{-1})+t(-2\theta_1,1), t \in \mathbb {R}\right\}.
$$
This set can be alternatively described as the set of $(x,y) \in M_0$ such that
$$
y=-\frac {1}{2\theta_1^2}+\frac {x}{2\theta_1}.
$$
Next we prove that the family $\{B_{\boldsymbol {v}(\theta_1)}: \theta_1>0\}$ is a partition of $M_0$.
Indeed, setting $z=\theta_1^{-1}$
we obtain that  $(x,y) \in B_{\boldsymbol {v}(\theta_1)}$ if and only if
$$
z^2+xz+2y=0,
$$
whose unique positive solution is given by $z=\frac {\sqrt {x^2-8y}-x}{2}$
(recall that $y<0$ when $(x,y)\in M_0$), and thus
$$
\theta_1=\frac {2}{\sqrt {x^2-8y}-x},
$$
and clearly, by varying $\theta_1>0$, a partition of $M_0$ is obtained.
As a consequence, setting
$$
\varphi_{\Theta}(x,y)=\boldsymbol {v}\left(\frac {2}{\sqrt {x^2-8y}-x}\right),\,\,\,\, (x,y) \in M_0
$$
we can apply Theorem \ref{eproj}. Since $A_{\boldsymbol {v}(\theta_1)}(\Theta)$ is the tangent direction 
to $\Theta$ at $\boldsymbol {v}(\theta_1)$, and $\Theta$ is the graph of a strictly convex function, it 
is not possible to find any positive $\theta_0\neq \theta_1$, such that $\boldsymbol {v}(\theta_0)$ lies 
in $\boldsymbol {v}(\theta_1)+A_{\boldsymbol {v}(\theta_1)}(\Theta)$. Except at $\theta_0$, the rate 
function $I (\boldsymbol {v}(\theta_{1}))$ in \eqref{cracontr} is always strictly smaller than 
$D(\boldsymbol {v}(\theta_{1})||\boldsymbol {v}(\theta_{0}))$, for any $\theta_{1}>0$.

Now let us take as a new subfamily the convex set
$$
\Theta_c=\left \{(\theta_1,\theta_2): \theta_2\geq \theta_1^2\right \},
$$
which corresponds to the assumption that $\sigma \geq |\mu|$. Its image under $\nabla \kappa$ is the set
$$
M_c=\left \{(\mu_1,\mu_2): \mu_2\leq -\mu_1^2\right \}.
$$
As long as $\boldsymbol {\mu} \in M_c$, the maximizer of $\ell (\boldsymbol {\theta}; \boldsymbol {\mu})$ in $\boldsymbol {\theta} \in \Theta_c$ is just the unconstrained one, that is $\boldsymbol {\theta}^{\ast}=\nabla \kappa^{\ast} (\boldsymbol {\mu})\in\Theta_c$. Now let $\boldsymbol {\mu}=(x,y) \in M_0 \setminus M_c$. If $x>0$, then it is easily verified that the maximizer lies on $\Theta$, that is on the boundary of $\Theta_c$, hence it is given again by
$$\boldsymbol {\theta}^{\ast}=\boldsymbol {v} \left( \frac {2}{\sqrt {x^2-8y}-x}\right).$$
By symmetry, for $x<0$, the maximizer is 
$$\boldsymbol {\theta}^{\ast}=\boldsymbol {v} \left( -\frac {2}{\sqrt {x^2-8y}+x}\right).
$$
By putting the three cases above together the continuous MLE $\varphi_{\Theta_c}$ is obtained.
From Theorem \ref{pyta1} we deduce the LDP for the laws of $\{\varphi_{\Theta_c} ({\bar {\boldsymbol {x}}_n}), n\in \mathbb{N}\}$ with rate function $D(\boldsymbol {\theta}||\boldsymbol {\theta}_0)$, for any $\boldsymbol {\theta} \in \Theta_{c}$, when the samples are drawn from $P_{\boldsymbol {\theta}_0}$, with $\boldsymbol {\theta}_0 \in \Theta_c$. This agrees with the fact that, for any $\theta_1>0$, the set of $\boldsymbol{\mu}$ such that the maximizer of 
$\ell (\boldsymbol {\theta}; \boldsymbol {\mu})$ in $\boldsymbol {\theta} \in \Theta_c$ is equal to $\boldsymbol {v}(\theta_1)$, is
$$\nabla\kappa(T_{\Theta_c}(\boldsymbol {v}(\theta_1)))=B_{\boldsymbol {v}(\theta_1)}\setminus int(M_c).$$
This is an expression of the fact that the  points of the line $B_{\boldsymbol {v}(\theta_1)}$ intersecting the interior of $M_c$ have to be removed from the set $\nabla\kappa(T_{\Theta_c}(\boldsymbol {v}(\theta_1)))$, since they are feasible mean values under the submodel with $\boldsymbol {\theta}\in\Theta_c$.

\end{example}

\section{Dual families and inverse LDP's}\label{sec:Letac}
This section is devoted to a rather more specific topic, for which we need to revise our notation: the cumulant generating function $\kappa$ of a measure $\lambda$ is denoted by $\kappa_{\lambda}$ in this section. Here is the key definition, which is slightly different from that given in Letac (2022), Section 3.1, inspired in turn by Barndorff-Nielsen (1978, Section 9.1, page 142).

\begin{definition} Let $\lambda$ be a $\sigma$-finite Borel measure on $\mathbb{R}^d$ with a steep cumulant generating function $\kappa_{\lambda}$. A $\sigma$-finite Borel measure $\lambda^{\ast}$ on $\mathbb{R}^d$, with 
a steep cumulant generating function $\kappa_{\lambda^{\ast}}$ as well, is called the dual measure of $\lambda$ if
$$\kappa_{\lambda^{\ast}}=\kappa_{\lambda}^{\ast},$$
where $\kappa_{\lambda}^{\ast}$ is the conjugate function of $\kappa_{\lambda}$, defined in \eqref{supachie}.	
\end{definition}

We warn the reader that not all $\sigma$-finite Borel measures on $\mathbb {R}^d$ have dual measures (see Letac (2022) for specific examples). The notion of dual measure has immediate consequences that we collect here:

\noindent
i) the dual measure $\lambda^{\ast}$ has the measure $\lambda$ as its dual; 

\noindent
ii) If $dom(\kappa_{\lambda})$ and $dom(\kappa_{\lambda^{\ast}})$ are the effective domains of $\lambda$ and $\lambda^{\ast}$ and $C(\lambda)$ and $C(\lambda^{\ast})$ are their convex supports, respectively, then:
$$
int(dom(\kappa_{\lambda^{\ast}}))=int(C(\lambda)), \,\,\, int(C(\lambda^{\ast}))=int(dom(\kappa_{\lambda}));
$$

\noindent
iii) $\nabla \kappa_{\lambda^{\ast}}$ is the inverse of $\nabla \kappa_{\lambda}$, from $int(C(\lambda))$ onto $dom(\kappa_{\lambda})$.

As a matter of fact, Letac (2022) gave directly property iii) as the definition of dual measure. In this way an additional affine term is allowed in its cumulant generating function $\kappa_{{\lambda}^{\ast}}$ which however leaves unchanged the natural exponential families generated by $\lambda$ and $\lambda^{\ast}$, to which we turn our attention.

When $\lambda$ has a dual measure there is an interesting connection between the LDP's concerning the full exponential families generated by $\lambda$ and $\lambda^{\ast}$, which we derive from the relation 
concening KL divergences within the exponential families 
$\{P_{\boldsymbol {\theta}},\boldsymbol \theta \in dom(\kappa_{\lambda})\}$ and 
$\{Q_{\boldsymbol {\mu}},\boldsymbol {\mu} \in  dom(\kappa_{\lambda^{\ast}})\}$, defined by
\begin{equation}\label{duals}
\frac {dP_{\boldsymbol {\theta}}}{d\lambda}(\boldsymbol x)=\exp\{\boldsymbol {\theta}\cdot \boldsymbol {x}-\kappa_{\lambda}(\boldsymbol {\theta})\}\ \mbox{and}\ \frac {dQ_{\boldsymbol {\mu}}}{d\lambda^{\ast}}(\boldsymbol x)=\exp\{\boldsymbol {\mu}\cdot \boldsymbol {x}-\kappa_{{\lambda}^{\ast}}(\boldsymbol {\mu})\},
\end{equation}
respectively. In what follows the KL divergence $D(Q_{\boldsymbol {\mu}_1}||Q_{\boldsymbol {\mu}_2})$ will be rewritten as $\tilde{D}(\boldsymbol {\mu}_1||\boldsymbol {\mu}_2)$ for any pair $\boldsymbol {\mu}_1,\boldsymbol {\mu}_2\in dom(\kappa_{{\lambda}^{\ast}})$.

\begin{proposition}\label{dualite} Let $\lambda$ and $\lambda^{\ast}$ be $\sigma$-finite Borel measures on $\mathbb {R}^d$ that are dual.
For $\boldsymbol {\theta}, \boldsymbol {\theta}_0 \in int(dom(\kappa_{\lambda}))$, define 
$$\boldsymbol {\mu}=\nabla \kappa_{\lambda} (\boldsymbol {\theta}),\,\,\,\, \boldsymbol {\mu}_0=\nabla \kappa_{\lambda} (\boldsymbol {\theta}_0) \in int(dom(\kappa_{{\lambda}^{\ast}})).$$
Then
$$D(\boldsymbol {\theta}||\boldsymbol {\theta}_0)=\tilde{D}(\boldsymbol {\mu}_0||\boldsymbol {\mu}).$$
\end{proposition}
\begin{proof}
The desired equality can be obtained from \eqref{againdiv} and \eqref{againdiv2}, since 
$$
D(\boldsymbol{\theta}||\boldsymbol{\theta}_0)
=\kappa_{\lambda} (\boldsymbol {\theta}_0)-\ell^{\ast}(\nabla \kappa_{\lambda} (\boldsymbol {\theta}); \boldsymbol {\theta}_0)=
\kappa_{{\lambda}^{\ast}}^{\ast}(\boldsymbol {\theta}_0)-\ell^{\ast}(\boldsymbol {\mu}; \boldsymbol {\theta}_0)=\tilde{D}(\boldsymbol {\mu}_0||\boldsymbol {\mu}),
$$
recalling that $\boldsymbol {\theta}_0=\nabla \kappa_{\lambda^{\ast}} (\boldsymbol {\mu}_0)$.
\end{proof}

From the above proposition the following corollary can be easily obtained. In order to avoid confusion, we continue to consider posterior distributions and MLE's for the natural parameters of both the dual exponential families defined in \eqref{duals}.

\begin{corollary}\label{corollary}
Suppose $\lambda$ is a regular $\sigma$-finite measure in $\mathbb{R}^d$
with $dom(\kappa_{\lambda})$ open. Let $\nu$ be a probability measure on $dom(\kappa_{\lambda})$, with a full relative support. Let $\{\bar{\boldsymbol{x}}_{n}\}$ be a sequence in $\mathbb {R}^d$ converging to $\boldsymbol {\mu}_0 \in int(C(\lambda))$, and let $\boldsymbol {\theta}_0$ be such that $\nabla \kappa_{\lambda} (\boldsymbol {\theta}_0)=\boldsymbol {\mu}_0$.

Next suppose that $\lambda$ has a dual measure $\lambda^{\ast}$. Consider i.i.d. random variables $\boldsymbol{x}^{\ast}_1,\ldots, \boldsymbol{x}^{\ast}_n$ drawn from $Q_{\boldsymbol {\mu}_0}$, defined in \eqref{duals}, and suppose that their sample mean $\bar{\boldsymbol{x}}_{n}^{\ast}$ takes values in $int(C(\lambda^*))=dom(\kappa_\lambda)$ with probability $1$. 

Then the sequence of posterior laws $\{\pi_n(\cdot|\bar{\boldsymbol{x}}_{n})\}$ on the natural parameter $\boldsymbol {\theta}$ defined in \eqref{bayes} (for the family $\{P_{\boldsymbol {\theta}}, \boldsymbol {\theta} \in \Theta_0\}$) and the sequence of laws of $\{\nabla \kappa_{\lambda}(\bar{\boldsymbol{x}}_{n}^{\ast})\}$ (the MLE
of the natural parameter $\boldsymbol {\mu}$ for the family $\{Q_{\boldsymbol {\mu}}, \boldsymbol {\mu} \in M_0\}$), have LDP's with rate functions $J$ and $\tilde I$ on $dom(\kappa_{\lambda})$ and $int(dom(\kappa_{{\lambda}^{\ast}}))$, respectively, related by
$$J(\boldsymbol {\theta})=D(\boldsymbol {\theta}_0||\boldsymbol {\theta})=\tilde{D}(\boldsymbol {\mu}||\boldsymbol {\mu}_0)=\tilde I(\boldsymbol {\mu}),$$
where $\boldsymbol \mu=\nabla \kappa_{\lambda} (\boldsymbol {\theta})$.
\end{corollary}

Differently from the choice for the presentation of the previous result, in the following classical example of dual measures we use the change of variable  $\boldsymbol \mu=\nabla \kappa_{\lambda} (\boldsymbol {\theta})$ in the posterior distribution, achieving the identity of the rate functions for the two LDP's.

\begin{example}\rm
The probability measure $\lambda$ equal to the Poisson law with mean $1$ generates the Poisson exponential family, and it is clearly regular. Its cumulant generating function is given by $\kappa_{\lambda}(\theta)=e^{\theta}-1$ and once a prior distribution $\omega$ is placed on the mean value parameter $\mu=e^{\theta}$ the distributions on $\mu$ conditional to $\bar {x}_n$ is given by 
$$
\omega_n(B|\bar {x}_n) =\frac  {\int_{B}e^{n({\bar {x}_n \cdot \log \mu}-\mu)}\omega(d\mu)}
{\int_0^{\infty}e^{n({\bar {x}_n \cdot \log \mu}-\mu)}\omega(d\mu) },
$$
where $B$ is any Borel subset of the positive real line. The function $\kappa_{\lambda}$ has the convex conjugate 
$$
\kappa^{\ast}_{\lambda}(\mu)=\left\{\begin{array}{ll}
\mu \log(\mu)-\mu+1&\ \mbox{for}\ \mu\geq 0\\
+\infty&\ \mbox{otherwise},
\end{array}\right.
$$ 
which is, up to the additive constant $1$, the cumulant generating function of the probability measure $\lambda^{\ast}$ with density
$$
f(y)=\frac {\pi}{4}\int_0^{+\infty} e^{\pi v/4}\cos(-v-vy+v\log v)dv,\,\,\, y \in \mathbb{R}
$$ 
with respect to the Lebesgue measure on the real line. This is the density of $-X-1$, where $X$ has the Landau distribution (see Landau 
(1967); see also Eaton et al. (1971) for informations about this law).

The family $\{Q_{\mu}, \mu>0\}$ appearing in the previous corollary, is the natural exponential family generated by $\lambda^{\ast}$. Its normalized 
log-likelihood function, as a function of the natural parameter $\mu$, is given by 
$$
\ell ({\mu}; \bar {x}_{n}^{\ast})=\mu \cdot \bar {x}_{n}^{\ast}-\kappa^{\ast}_{\lambda}(\mu),
$$
which is maximized by
$$\mu=\nabla\kappa_\lambda(\bar {x}_{n}^{\ast})=e^{\bar {x}_{n}^{\ast}}.$$
As a result, for $x^{\ast}_1,\ldots, x^{\ast}_n$ i.i.d. from $Q_{\mu_0}$, as $n \to \infty$ the sequence of laws of the MLE
$\{e^{\bar {x}_{n}^{\ast}}\}$ obeys the same LDP as the sequence of posterior laws $\{\omega_n(\cdot|\bar {x}_n)\}$ on the positive real line, 
when $\{\bar {x}_n\}$ is a sequence converging to $\mu _0>0$. For any $\mu_{0}>0$, the rate function for the MLE of the parameter $\mu$ is
$$
\tilde I(\mu)=\tilde{D}(\mu||\mu_0)=(\mu-{\mu}_0)\nabla \kappa_{\lambda^{\ast}} (\mu)+\kappa_{\lambda^{\ast}}(\mu)-\kappa_{\lambda^{\ast}}(\mu_0)=
\mu_0\log(\mu_0/{\mu})+\mu-\mu_0,\quad \mu>0.
$$
By Theorem \ref{invsanofi} the rate function $J$ for the posterior laws on the parameter $\theta$ is obtained as
$$J(\theta)=D(\theta_0||\theta)=(\theta_0-\theta)e^{\theta_0}-e^{\theta_0}+e^{\theta},\quad \theta\in\mathbb{R}$$
where $e^{\theta_0}=\mu_0$. A trivial application of the contraction principle allows to deduce that the rate function for $\{\omega_n(\cdot|\bar {x}_n)\}$ is
$$J(\log\mu)=\tilde I(\mu),\quad \mu>0.$$
\end{example}

\section{LDP's for posteriors under more general conditions}\label{sec:dom(K)-not-open}

The aim of this section is to explain the circumstances under which the LDP stated in Theorem \ref{invsanofi} continues to hold when 
the essential domain $dom(\kappa)$ of the cumulant generating function $\kappa$ of the reference measure $\lambda$ is not open.
In this case $\kappa$ remains continuous in the interior of $dom(\kappa)$, but this is not necessarily true 
at boundary points of $dom(\kappa)$. To deal with these situations we introduce the following notion.
\begin{definition}
Let $B \subset dom(\kappa)$ and $\boldsymbol {\theta} \in B\cap \partial dom(\kappa)$. We say that $\boldsymbol {\theta}$ \emph{is a continuity point for} $\kappa$ on $B$ if for any sequence $\{\boldsymbol{\theta}_{\ell}\} \subset B\cap int(dom (\kappa))$ such that $\boldsymbol{\theta}_{\ell} \to \boldsymbol{\theta}$ it happens that $\kappa (\boldsymbol{\theta}_{\ell}) \to \kappa(\boldsymbol{\theta})$ as $\ell \to \infty$.
\end{definition}

Additionally we want to allow the state space of the LDP to be what we find convenient to call a \emph{measurable support} of $\lambda$,
according to the following definition. Recall that for any probability measure $\nu$ on $dom(\kappa)$, $S(\nu)$ is the topological support of $\nu$ in the relative topology of $dom(\kappa)$.

\begin{definition}
Let $\nu$ be a probability measure $\nu$ such that $\nu(dom (\kappa))=1$. If $\Theta\subset S(\nu)$ and $\nu(\Theta)=1$, we say that $\Theta$ is a measurable support of $\nu$.
\end{definition}

\begin{theorem}\label{cor2} Let $\lambda$ be any $\sigma$-finite Borel measure on $\mathbb{R}^d$ which is not concentrated on a proper affine submanifold and whose cumulant generating function $\kappa$ has an essential domain $dom(\kappa)$ with non-empty interior. Moreover let $\nu$ be a probability measure on $dom(\kappa)$ and $\Theta\subset S(\nu)$ be a measurable support of $\nu$. Assume that the following conditions hold.
\begin{itemize}
\item a. The sequence  $\{\bar{\boldsymbol{x}}_{n}, n\in \mathbb{N}\}$ converges to $\boldsymbol{\mu}_0 \in int(C(\lambda))$.
\item b. There exists an $m$-projection $\Pi_{S(\nu)}(\boldsymbol {\theta}_0)$ of $\boldsymbol {\theta}_{0}=\nabla\kappa^*(\boldsymbol {\mu}_{0})$ on $S(\nu)$ which is either in $int(dom(\kappa))$, or it is a continuity point for $\kappa$ on $S(\nu)$.
\item c. Any $\boldsymbol{\theta} \in \Theta \cap \partial dom(\kappa)$ is a continuity point for $\kappa$ on $\Theta$.
\end{itemize}
Then the sequence of probability measures $\{\pi_{n}(\cdot|\bar {\boldsymbol{x}}_{n}), n \in \mathbb {N}\}$ on $\Theta$, defined by \eqref{bayes}, satisfies a LDP with rate function
$$J(\boldsymbol{\theta})=D(\boldsymbol{\theta}_0||\boldsymbol{\theta})-D(\boldsymbol{\theta}_0||\Pi_{S(\nu)}(\boldsymbol {\theta}_0)),\,\,\, \boldsymbol{\theta} \in \Theta.$$
\end{theorem}

A few remarks before proving Theorem \ref{cor2}. Assumption a. remains unchanged with respect to Theorem \ref{invsanofi}: the sequence  $\{\bar{\boldsymbol{x}}_{n}, n\in \mathbb{N}\}$ has to converge to $\boldsymbol{\mu}_0 \in int(C(\lambda))$. Assumptions b. and c. depend on the prior probability $\nu$, but they always hold for one-dimensional exponential families, by lower semi-continuity of the function $\kappa$. The $m$-projection $\Pi_{S(\nu)}(\boldsymbol {\theta}_0)$ mentioned in b. exists by the fact that $S(\nu)$ is relatively closed, $\kappa$ is lower semi-continuous and has bounded (super)-level sets, but we need assumption b. to exclude that $\Pi_{S(\nu)}(\boldsymbol {\theta}_0)$ lies on the boundary of $S(\nu)$ without being a continuity point. Finally notice that it is not required that $\Pi_{S(\nu)}(\boldsymbol {\theta}_0)$ belongs to $\Theta$.

\begin{proof} As far as the upper bounds in Lemma \ref{th1} and Theorem \ref{invsanofi} are concerned the assumption about $\boldsymbol{\mu}_0$ continue to work as in the regular case. Indeed, for any $B \subset \mathbb{R}^d$, the convex function $\kappa^{\ast}_B$ is continuous in the interior of its essential domain, which necessarily includes the interior of $C(\lambda)$.

As far as the lower bounds in  Lemma \ref{th1} and Theorem \ref{invsanofi} are concerned, we have to take into account that  $\Pi_{S(\nu)}(\boldsymbol {\theta}_0)$ and $\boldsymbol {\theta}_{\ast}$ may lie in the boundary of $dom(\kappa)$. However, assumptions b. and c. ensure that for any $m \in \mathbb{N}$, one can replace $\Pi_{S(\nu)}(\boldsymbol {\theta}_0)$ and $\boldsymbol {\theta}_{\ast}$ with $\boldsymbol {\theta}_{\nu}^{m} \in S(\nu)$ and $\boldsymbol {\theta}_{\ast}^{m} \in T$ with the property
$$
\ell(\boldsymbol {\theta}^{m}_{\nu};\boldsymbol{\mu}_0)>\ell(\boldsymbol {\theta}_{\nu};\boldsymbol{\mu}_0)-1/m,\,\,\,
\ell(\boldsymbol {\theta}^{m}_{\ast};\boldsymbol{\mu}_0)>\ell(\boldsymbol {\theta}_{\ast};\boldsymbol{\mu}_0)-1/m,
$$
respectively. Now one repeats the proof of \eqref{lb1} and \eqref{lb2} to prove these lower bounds with the right hand side equal to $\ell(\Pi_{S(\nu)}(\boldsymbol {\theta}_0);\boldsymbol{\mu}_0)-1/m$ and $\ell(\boldsymbol {\theta}_{\ast};\boldsymbol{\mu}_0)-1/m$, respectively. Since they hold irrespectively of $m$, they continue to hold when $m^{-1}$ is replaced by zero.

\end{proof}

The following classical example (see e.g. Barndorff-Nielsen (1978), Example 7.3, page 104), serves as an illustration of the previous result.

\begin{example}\rm
Consider the probability measure $\mu$ on $\mathbb{R}^2$ given by
$$\lambda(dx_1,dx_2)=\frac {1}{2 \sqrt {\pi}(1+x_1^2)^{3/2}}\exp \left \{-\left(x_1^2+\frac {x_2^2}{4(1+x_1^2)}\right)\right \}dx_1dx_2.$$
The essential domain of its cumulant generating function $\kappa$ is given by
$$
dom (\kappa)= \{(\theta_1,\theta_2) \in \mathbb{R}^2:\theta_2 \in (-1,1)\} \cup \{(0, \pm1)\}.
$$
The peculiarity of this example is that on the boundary points $(0,\pm 1)$ of $dom(\kappa)$ the function $\kappa$ is finite but it is not continuous. Indeed on the curve  
\begin{equation}\label{curve}
\theta_2=\sqrt{1-\theta_1^3},\,\,\,\, \theta_1 \in (0,1]
\end{equation}
the function $\kappa$ tends to $+\infty$ as $\theta_1 \to 0$. Now define $\eta(u)=(u,\sqrt{1-u^3})$ for $u \in (0,1]$, let $\gamma$ be the uniform prior on $U=(0,1]$ and let $\nu=\gamma \circ \eta^{-1}$ be its image measure on $dom(\kappa)$. Clearly $\Theta=\eta(U)$ is a measurable support of $\nu$, whereas its topological support $S(\nu)$ contains also the limit point $(0,1)$. We are going to apply Theorem \ref{cor2} to this example. Since $\eta$ is a homeomorphism of $U$ onto $\Theta$, the LDP can be immediately transferred to $U$.

The measure $\lambda$ has the whole $\mathbb {R}^2$ as support, so assuming that $\bar {\boldsymbol{x}}_{n}$ converges to some $\boldsymbol{\mu}_0=({\mu}_{01},{\mu}_{02})\in \mathbb{R}^2$ guarantees assumption a. of the theorem. As a consequence the function $\ell(\eta(u);\boldsymbol{\mu}_0)$ defined for $u \in U$, which is continuous in any interval $[\varepsilon,1]$ with $\varepsilon >0$ and tends to $-\infty$ as $u \to 0$, has a maximizer in $u$, say $\tilde u$, which depends on $\boldsymbol{\mu}_0$, whose image under $\eta (\tilde u)$ is an interior point of $dom(\kappa)$. So assumption b. is automatic. As far as assumption c. is concerned nothing has to be verified since $M$ does not contain any boundary point of $dom(\kappa)$. As a consequence Theorem \ref{cor2} can be applied, with the rate function given by
\begin{multline*}
\iota_U(u)=\ell(\eta(\tilde u);\boldsymbol{\mu}_0)-\ell(\eta(u);\boldsymbol{\mu}_0)\\
=(\tilde u-u) \mu_{0,1}+(\sqrt{1-\tilde u^3}-\sqrt{1-u^3}) \mu_{0,2}-\kappa (\tilde u, \sqrt{1-\tilde u^3})+\kappa (u, \sqrt{1-u^3}),
\end{multline*}
where $u \in M$. In particular the rate function tends to $+ \infty$ as $u \to 0$: indeed the rate function does not see at all that the cumulant generating function is finite at the boundary point $(0,1)$, since this is not contained in the image of $U$ under $\eta$. 

Now suppose that $\Theta=\eta(\overline U)$, hence $u=0$ is added and mapped into the limit point $(0,1)$. Then Theorem \ref{cor2} cannot be applied because of the failure of assumption c. As a matter of fact the LDP fails on spheres of radius $\varepsilon>0$ suitably small around $(0,1)$, whose posterior probability approaches zero with an exponential rate that can be made arbitrarily large by chosing $\varepsilon$ small enough, which is not compatible with the finite value $\iota_{\bar U}(0)=\ell((0,1);\boldsymbol{\mu}_0)$.

Finally notice that if the curve \eqref{curve} is replaced by
$$\theta_2=1-\theta_1,\,\,\,\, \theta_1 \in (0,1]$$
again with a uniform prior on $\theta_1$, the above phenomenon disappears, due to the continuity of $\kappa$ on the curve, up to the boundary point $(0,1)$, inherited by lower semi-continuity of $\kappa$ on the whole of $\mathbb {R}^2$.
\end{example}

\paragraph{Acknowledgements.}
The authors thank the anonymous referees for many useful suggestions. They also thank Steffen Lauritzen for relevant suggestions about convex subfamilies of exponential families and for the reference to Chentsov (1966). Finally, they acknowledge the generosity of G\'{e}rard Letac for having provided them an earlier draft of Letac (2022).


\begin{thebibliography}{spc}
\bibitem{Arcones}
Arcones, M.A., 2006. Large deviations for M-estimators. Ann. Inst. Statist. Math. 58, 21--52.
\bibitem{BahadurZabellGupta}
Bahadur, R.R., Zabell, S.L., Gupta, J.C., 1980. Large deviations, tests, and estimates. 
Asymptotic theory of statistical tests and estimation (Proc. Adv. Internat. Sympos., 
Univ. North Carolina, Chapel Hill, N.C., 1979), pp. 33--64, Academic Press, 
New York-London-Toronto.
\bibitem{Barndorff-Nielsen}
Barndorff-Nielsen, O., 1978. Information and exponential families
in statistical theory. John Wiley \& Sons, Ltd., Chichester.
\bibitem{Brown}
Brown, L.D., 1986. Fundamentals of statistical exponential families
with applications in statistical decision theory. Institute of 
Mathematical Statistics, Hayward.
\bibitem{Chentsov}
Chentsov, N.N., 1966. A systematic theory of exponential families 
of probability distributions. Teor. Verojatnost. i Primenen 11, 483--494
\bibitem{Csiszar}
Csisz\'ar, I., 1975. $I$-divergence geometry of probability distributions 
and minimization problems. Ann. Probability 3, 146--158. 
\bibitem{DemboZeitouni}
Dembo, A., Zeitouni, O., 1998. Large deviations techniques and 
applications. Second edition. Springer-Verlag, New York.
\bibitem{EatonMorrisRubin}
Eaton, M.L., Morris, C., Rubin, H., 1971. On extreme stable laws and
some applications. J. Appl. Probability 8, 794--801.
\bibitem{GaneshOConnell1999}
Ganesh, A., O'Connell, N., 1999. An inverse of Sanov's theorem. Statist.
Probab. Lett. 42, 201--206.
\bibitem{GaneshOConnell2000}
Ganesh, A.J., O'Connell, N., 2000. A large-deviation principle for 
Dirichlet posteriors. Bernoulli 6, 1021--1034.
\bibitem{GhosalGhosHRamamoorthi}
Ghosal, S., Ghosh, J.K., Ramamoorthi, R.V., 1999. Consistency issues in 
Bayesian nonparametrics. Asymptotics, nonparametrics, and time series, 
pp. 639--667, Statist. Textbooks Monogr., 158, Dekker, New York.
\bibitem{KesterKallenberg}
Kester, A.D.M., Kallenberg, W.C.M., 1986. Large deviations of estimators.
Ann. Statist. 14, 648--664.
\bibitem{KleijnvanderVaart}
Kleijn, B.J.K., van der Vaart, A.W., 2006. Misspecification in 
infinite-dimensional Bayesian statistics. Ann. Statist. 34, 837--877.
\bibitem{Laundau}
Landau, L.D., 1967. Collected papers of L.D. Landau. Edited and with an 
introduction by D. ter Haar. Second printing Gordon and Breach Science 
Publishers, New York-London-Paris.
\bibitem{Letac}
Letac, G., 2022. Duality for real and multivariate exponential families.
J. Multivariate Anal. 188, Paper No. 104811, 22 pp.
\bibitem{Macci}
Macci, C., 2014. Extension of some large deviation results for posterior 
distributions. J. Korean Statist. Soc. 43, 189--200.
\bibitem{MacciPetrella}
Macci, C., Petrella, L., 2009. Censored exponential data: large deviations
for MLEs and posterior distributions. Comm. Statist. Theory Methods 38, 
2435--2452.
\bibitem{MiaoHahn}
Miao, W., Hahn, M., 1997. Existence of maximum likelihood estimates for multi-dimensional exponential families. Scand. J. Statist. 24, 371--386.
\bibitem{Morris}
Morris, C.N., 1982. Natural exponential families with quadratic variance functions. Ann. Statist. 10, 65--80. 
\bibitem{Nielsen}
Nielsen, F., 2018. What is ... an information projection?
Notices Amer. Math. Soc. 65, 321--324.
\bibitem{RobertsVarberg}
Roberts, A.W., Varberg, D.E., 1973. Convex functions. Academic Press, 
New York-London, 1973.
\bibitem{Simon}
Simon, G., 1973. Additivity of information in exponential family probability
laws. J. Amer. Statist. Assoc. 68, 478--482.
\bibitem{Wasserman}
Wasserman, L., 1998.
Asymptotic properties of nonparametric Bayesian procedures. Practical 
nonparametric and semiparametric Bayesian statistics, pp. 293--304,
Lect. Notes Stat., 133, Springer, New York.
\end{thebibliography}
\end{document}